\documentclass[aop,MSNbibl,number,citesort,dvips]{arximspdf}

\doi{10.1214/11-AOP715}
\volume{41}
\issue{3B}
\pubyear{2013}
\firstpage{2182}
\lastpage{2224}

\makeatletter
\newcommand{\skp}[2]{\langle #1 , #2 \rangle}
\newcommand{\mv}[1]{\langle #1 \rangle}
\newcommand{\R}{\mathbb{R}}
\newcommand{\N}{\mathbb{N}}
\newcommand{\h}{{\fraca{1}{2}}}

\newtheorem{theorem}{Theorem}[section]
\newproclaim{definition}[theorem]{Definition}
\newtheorem{proposition}[theorem]{Proposition}
\newtheorem{lemma}[theorem]{Lemma}
\newproclaim{remark}[theorem]{Remark}

\newcommand{\id}{\mathrm{id}}
\newcommand{\Real}{\operatorname{Re}}
\newcommand{\Hess}{\operatorname{Hess}}
\newcommand{\osc}{\operatorname{osc}}
\newcommand{\var}{\operatorname{var}}
\newcommand{\cov}{\operatorname{cov}}
\newcommand{\im}{\operatorname{im}}
\newcommand{\eqref}[1]{(\ref{#1})}
\newcommand{\overset}{\stackrel}
\newcommand{\mathclap}{\scriptsize}
\newcommand{\fracd}[2]{({#1}/{#2})}
\newcommand{\fracc}[2]{{#1}/{(#2)}}
\newcommand{\fraca}[2]{{#1}/{#2}}
\makeatother

\begin{document}
\begin{frontmatter}

\title{Uniform logarithmic Sobolev inequalities for conservative spin
systems with super-quadratic single-site potential}
\runtitle{LSI for conservative spin systems}

\begin{aug}
\author[A]{\fnms{Georg} \snm{Menz}\corref{}\ead[label=e1]{Georg.Menz@mis.mpg.de}\thanksref{tt1}}
\and
\author[A]{\fnms{Felix} \snm{Otto}\ead[label=e2]{Felix.Otto@mis.mpg.de}\ead[label=u1,url]{http://www.mis.mpg.de/}}
\runauthor{G. Menz and F. Otto}
\affiliation{Max Planck Institute for Mathematics in the Sciences}
\address[A]{Inselstra{\ss}e 22\\ 04103 Leipzig \\ Germany \\
\printead{e1} \\
\hphantom{E-mail: }\printead*{e2}\\ \printead{u1}} %adresu isvedimo
%komanda gale!
\end{aug}
\thankstext{tt1}{Supported through the Gottfried Wilhelm Leibniz program,
the Bonn International Graduate School in Mathematics and the Max
Planck Institute for Mathematics in the Sciences in Leipzig.}

% HISTORY:
\received{\smonth{3} \syear{2011}}
\revised{\smonth{9} \syear{2011}}

% ABSTRACT
%
\begin{abstract}
We consider a noninteracting unbounded spin system with conservation
of the mean spin. We derive a uniform logarithmic Sobolev inequality
(LSI) provided the single-site potential is a bounded perturbation of a
strictly convex function. The scaling of the LSI constant is optimal in
the system size. The argument adapts the two-scale approach of
Grunewald, Villani, Westdickenberg and the second author from the
quadratic to the general case. Using an asymmetric Brascamp--Lieb-type
inequality for covariances, we reduce the task of deriving a uniform
LSI to the convexification of the coarse-grained Hamiltonian, which
follows from a general local Cram\'er theorem.
\end{abstract}

% KEYWORDS
%

\end{frontmatter}

%s1 ###
%s1 #&#
\section{Introduction and main result} \label{s_main_result}

The grand canonical ensemble $\mu$ is a probability measure on $\R^N$
given by
\[
% \label{e_cannonical_ensemble}
\mu(dx):= \frac{1}{Z} \exp ( -H(x)  )\,dx.
\]
Throughout the article, $Z$ denotes a generic normalization constant.
The value of $Z$ may change from line to line or even within a line.
The noninteracting Hamiltonian $H\dvtx  \R^N \to\R$ is given by a sum of
single-site potentials $\psi\dvtx  \R\to \R$ that are specified later,
that is,
%
%e1 ###
%e1 #&#
\begin{equation} \label{e_Hamiltonian}
H ( x  ) := \sum_{i=1}^N \psi(x_i).
\end{equation}
For a real number $m$, we consider the $N-1$ dimensional hyper-plane
$X_{N,m}$ given by
\[
% \label{e_definition_of_X_N_M}
X_{N,m} := \Biggl \{ x \in\R^N,  \frac{1}{N} \sum_{i=1}^N x_i =m
 \Biggr\}.
\]
We equip $X_{N,m}$ with the standard scalar product induced by $\R^N$, namely
\[
\skp{x}{\tilde x} := \sum_{i=1}^N x_i \tilde x_i .
\]
The restriction of $\mu$ to $X_{N,m}$ is called canonical ensemble
$\mu_{N,m}$, that is,
%
%e2 ###
%e2 #&#
\begin{equation}
\label{e_definition_of_mu_N_m}
\mu_{N,m}(dx) := \frac{1}{Z}  \exp ( - H(x)  )
\mathcal{H}_{ \lfloor X_{N,m} }^{N-1} (dx). %_{\lfloor \{
\end{equation}
Here, $\mathcal{H}_{\lfloor X_{N,m} }^{N-1}$ denotes the $N-1$
dimensional Hausdorff measure restricted to the hyperplane $X_{N,m}$.
For convenience, we introduce the notation:
\begin{eqnarray*}
a &\lesssim& b \quad\Leftrightarrow \quad  \mbox{there is a constant $C >
0 $ uniformly in the systems size $N$ and}\\[-5pt]
&&\hphantom{b \quad\Leftrightarrow \quad\;} \mbox{the mean spin $m$ such that } a \leq C b;\\
 a &\sim& b \quad \Leftrightarrow  \quad  \mbox{it holds that } a \lesssim b
\mbox{ and } b \lesssim a.
\end{eqnarray*}

In 1993, Varadhan (\cite{Var}, Lemma~5.3 ff.) posed the question for
which kind of single-site potential $\psi$ the canonical ensemble $\mu
_{N,m}$ satisfies a spectral gap inequality (SG) uniformly in the
system size $N$ and the mean spin $m$. A partial answer was given by
Caputo \cite{Cap}:
%
%th1.1 #&#
\begin{theorem}[(Caputo)]\label{p_caputo}
Assume that for the single-site potential $\psi$ there exist a
splitting $\psi= \psi_c + \delta\psi$ and constants $\beta_{-}$,
$\beta_{+} \in [ 0, \infty )$ such that for all $x \in
 [ 0, \infty ),$
%
%e3 ###
%e3 #&#
\begin{eqnarray}
\label{e_cond_cap}
&&\psi_c ''( x) \sim|x|^{\beta_{+}} + 1, \qquad  \psi_c''(- x) \sim
|x|^{\beta_{-}} + 1 \quad \mbox{and}\nonumber
\\[-8pt]
\\[-8pt]
&&|\delta\psi|+  |
\delta\psi'  | +  | \delta\psi''  | \lesssim1.\nonumber
\end{eqnarray}
Then the canonical ensemble $\mu_{N,m}$ satisfies the SG with constant
$\varrho>0$ uniformly in the system size $N$ and the mean spin $m$.
More precisely, for any function $f$,
\[
\var_{\mu_{N,m}}(f) = \int \biggl( f - \int f \,d \mu_{N,m}
\biggr)^2\,d\mu_{N,m} \leq\frac{1}{\varrho} \int|\nabla f|^2 \,d \mu_{N,m}.
\]
Here, $\nabla$ denotes the gradient determined by the Euclidean
structure of $X_{N,m}$.
\end{theorem}

In this article, we give a full answer to the question by
Varadhan \cite{Var} and also show that the last theorem can be
strengthened to the logarithmic Sobolev inequality (LSI).
%
%de1.2 #&#
\begin{definition}[(LSI)] \label{d_LSI}
Let $X$ be a Euclidean space. A Borel probability measure $\mu$ on $X$
satisfies the LSI with constant $\varrho>0$, if for all functions $f
\geq0$
%
%e4 ###
%e4 #&#
\begin{equation}\label{e_definition_of_LSI}
\int f \log f \,d \mu- \int f d\mu \log \biggl( \int f d\mu \biggr)
\leq\frac{1}{2 \varrho} \int\frac{|\nabla f|^2}{f}\,d\mu.
\end{equation}
Here, $\nabla$ denotes the gradient determined by the Euclidean
structure of $X$.
\end{definition}

%
%re1.3 #&#
\begin{remark}[(Gradient on $X_{N,m}$)]
If we choose $X=X_{N,m}$ in Definition~\ref{d_LSI}, we can calculate
$|\nabla f|^2$ in the following way: Extend $f\dvtx X_{N,m} \to\R$ to be
constant on the direction normal to $X_{N,m}$. Then %$f:\R^N \to\R$ and
\[
|\nabla f|^2 = \sum_{i=1}^N  \biggl| \frac{d}{dx_i} f  \biggr|^2.
\]
\end{remark}

The LSI was originally introduced by Gross \cite{Gross}. It yields the
SG and can be used as a powerful tool for studying spin systems. Like
the SG, the LSI implies exponential convergence to equilibrium of the
naturally associated conservative diffusion process. The rate of
convergence is given by the LSI constant $\varrho$; cf. \cite{R}, Chapter~3.2, and Remark~\ref{r_G_K}. Therefore, an appropriate
scaling of the LSI constant in the system size indicates the absence of
phase transitions. The SG yields convergence in the sense of variances
in contrast to the LSI, which yields convergence in the sense of
relative entropies. The SG and the LSI are also useful for deducing the
hydrodynamic limit; see \cite{Var} for the SG and \cite{GORV} for the LSI.

We consider three cases of different potentials: sub-quadratic,
quadratic and super-quadratic single-site potentials. In the case of
sub-quadratic single-site potentials, Barthe and Wolff \cite{BW} gave
a counterexample where the scaling in the system size of the SG and the
LSI constant of the canonical ensemble differs in the system size. More
precisely, they showed:
%
%th1.4 #&#
\begin{theorem}[(Barthe and Wolff)]
Assume that the single-site potential $\psi$ is given by
\[
\psi(x)=
\cases{
x ,&\quad  for   $x >0$,\cr
\infty ,&\quad  else.
}
\]
Then the SG constant $\varrho_1$ and the LSI constant $\varrho_2$ of
the canonical ensemble $\mu_{N,m}$ satisfy
\[
\varrho_1 \sim\frac{1}{m^2} \quad\mbox{and} \quad\varrho_2
\sim\frac{1}{Nm^2}.
\]
\end{theorem}

In the case of perturbed quadratic single-site potentials it is known
that Theorem~\ref{p_caputo} can be improved to the LSI. More
precisely, several authors (cf. \cite{LPY,Cha,GORV}) deduced the
following statement by different methods:\vadjust{\goodbreak}
%
%th1.5 #&#
\begin{theorem}[(Landim, Panizo and Yau)]\label{p_GORV}
Assume that the single-site potential $\psi$ is perturbed quadratic in
the following sense: There exists a splitting $\psi= \psi_c + \delta
\psi$ such that
%
%e5 ###
%e5 #&#
\begin{equation}
\label{e_cond_LPY}
\psi_c '' = 1 \quad\mbox{and} \quad|\delta\psi|+  | \delta
\psi'  | +  | \delta\psi''  | \lesssim1.
\end{equation}
Then the canonical ensemble $\mu_{N,m}$ satisfies the LSI with
constant $\varrho>0$ uniformly in the system size $N$ and the mean
spin $m$.
\end{theorem}

There is only left to consider the super-quadratic case. It is
conjectured that the optimal scaling LSI also holds if the single-site
potential $\psi$ is a bounded perturbation of a strictly convex
function; cf. \cite{LPY}, page~741, \cite{Cha}, Theorem~0.3~f., and
\cite{Cap}, page~226. Heuristically, this conjecture seems reasonable:
Because the LSI is closely linked to convexity (consider, e.g.,
the Bakry--\'Emery criterion), a perturbed strictly convex potential
should behave no worse than a perturbed quadratic one. However
technically, the methods for the quadratic case are not able to handle
the perturbed strictly convex case because they require an upper bound
on the second derivative of the Hamiltonian. In the main result of the
article we show that the conjecture from above is true:
%
%th1.6 #&#
\begin{theorem}\label{p_main_result}
Assume that the single-site potential $\psi$ is perturbed strict\-ly
convex in the sense that there is a splitting $\psi= \psi_c + \delta
\psi$ such that
%
%e6 ###
%e6 #&#
\begin{equation}
\label{e_perturbation_convex}
\psi_c'' \gtrsim1 \quad \mbox{and} \quad|\delta\psi| + |\delta
\psi' | \lesssim1.
\end{equation}
Then the canonical ensemble $\mu_{N,m}$ satisfies the LSI with
constant $\varrho>0$ uniformly in the system size $N$ and the mean
spin $m$.
\end{theorem}

%re1.7 #&#
\begin{remark}[(From Glauber to Kawasaki)] \label{r_G_K}
The bound on the r.h.s. of \eqref{e_definition_of_LSI} is given in
terms of the Glauber dynamics in the sense that we have endowed
$X_{N,m}$ with the standard Euclidean structure inherited from $\R^N$.
By the discrete Poincar\'e inequality, one can recover the bound for
the Kawasaki dynamics (cf. \cite{GORV}, Remark 15, or \cite{Cap}) in
the sense that one endows $X_{N,m}$ with the Euclidean structure coming
from the discrete $H^{-1}$-norm. More precisely, if $\Lambda$ is a
cubic lattice in any dimension of width $L$, then Theorem~\ref
{p_main_result} yields the LSI for Kawasaki dynamics with constant
$L^{-2} \varrho$, which is the optimal scaling in $L$; cf. \cite{Y2}.
\end{remark}

Note that the standard criteria for the SG and the LSI (cf. \hyperref[s_basic_facts_about_LSI]{Appendix}) fail for the canonical ensemble $\mu_{N,m}$:
\begin{itemize}
\item The \textit{Tensorization principle} for the SG and the LSI does not apply
because of the restriction to the hyper-plane $X_{N,m}$; cf. \cite{GZ}, Theorem~4.4, or Theorem~\ref{p_tensorization_principle}.
\item The \textit{Bakry--\'Emery} criterion does not apply
because the Hamiltonian $H$ is not strictly convex; cf. \cite{BE}, Proposition 3 and Corollary 2, or
Theorem~\ref{p_criterion_bakry_emery}.\vadjust{\goodbreak}
\item The \textit{Holley--Stroock} criterion does not help
because the LSI constant $\varrho$ has to be independent of the system
size $N$; cf. \cite{HS}, page~1184, or Theorem~\ref{p_criterion_holley_stroock}.
%because of the restriction to the hyper-plane $X_{N,m}$.
\end{itemize}
Therefore, a more elaborated machinery was needed for the proof of
Theorems~\ref{p_caputo}
and~\ref{p_GORV}. The approach of Caputo to Theorem~\ref{p_caputo}
seems to be restricted to the SG because it relies on the spectral
nature of the SG.
For the proof of Theorem~\ref{p_GORV}, Landim, Panizo and Yau \cite
{LPY} and
Chafa\"{\i} \cite{Cha} used the Lu--Yau martingale method that was
originally introduced in \cite{LY}
to deduce an analog version of Theorem~\ref{p_GORV} in the case of
discrete spin values. Recently,
Grunewald, Villani, Westdickenberg and the second author \cite{GORV}
provided a new technique for deducing Theorem~\ref{p_GORV},
called the two-scale approach. We follow this approach in the proof of
Theorem~\ref{p_main_result}.

The limiting factor for extending Theorem~\ref{p_GORV} to more general
single-site potentials is almost the same for the Lu--Yau martingale
method and for the two-scale approach: It is the estimation of a
covariance term w.r.t. the measure $\mu_{N,m}$ conditioned on a
special event; cf. \cite{LPY}, (4.6), and \cite{GORV}, (42). In the
two-scale approach one has to estimate for some large but fixed $K \gg
1$ and any nonnegative function $f$ the covariance
\[ %\label{e_critical_covariance_untransformed}
 \Biggl| \cov_{\mu_{K,m}}  \Biggl( f, \frac{1}{K} \sum_{i=1}^{K} \psi
'(x_i)  \Biggr)  \Biggr|.
\]
In \cite{GORV}, this term term was estimated by using a standard
estimate (cf. Lemma~\ref{p_cov_clas} and \cite{GORV}, Lemma~22) that
only can be applied for perturbed quadratic single-site potentials
$\psi$. We get around this difficulty by making the following
adaptations: Instead of one-time coarse-graining of big blocks, we
consider iterative coarse-graining of pairs. As a consequence we only
have to estimate the covariance term from above in the case $K=2$.
Because $\mu_{2,m}$ is a one-dimensional measure, we are able to apply
the more robust asymmetric Brascamp--Lieb inequality (cf. Lemma~\ref
{p_asymmetric_BL}) that can also be applied for perturbed strictly
convex single-site potentials $\psi$.

Recently, the optimal scaling LSI was established in \cite{Menz} by
the first author for a weakly interacting Hamiltonian with perturbed
quadratic single-site potentials $\psi$, that is,
\[
H(x) = \sum_{i=1}^N \psi(x_i) + \varepsilon\sum_{1\leq i<j \leq N}
b_{ij} x_i x_j.
\]
Because the original two-scale approach was used, it is an interesting
question if one could extend this result to perturbed strictly convex
single-site potentials. A~direct transfer of the argument of \cite
{Menz} fails because of the iterative structure of the proof of
Theorem~\ref{p_main_result}.% Furthermore, the argument of \cite{Menz}
%for the strict convexity of the coarse-grained Hamiltonian for small $

%

The remaining part of this article is organized as follows. In
Section~\ref{s_core_of_proof_of_main_result} we prove the main result.
The auxiliary results of Section~\ref{s_core_of_proof_of_main_result}\vadjust{\goodbreak}
are proved in Section~\ref{s_auxiliary_results}. There is one
exception: The convexification of the single-site potential by iterated
renormalization (see Theorem~\ref{p_convexification}) is proved in
Section~\ref{s_convexification}. In the short \hyperref[s_basic_facts_about_LSI]{Appendix} we state the standard criteria for the SG and
the LSI.\vspace*{-1pt}

\section{Adapted two-scale approach}\label{s_main_result_step_1}\vspace*{-1pt}
%s2.1 ###
%s2.1 #&#
\subsection{Proof of the main result} \label{s_core_of_proof_of_main_result}
The proof of Theorem~\ref{p_main_result} is based on an adaptation of
the two-scale approach of \cite{GORV}. We start with introducing the
concept of coarse-graining of pairs. We recommend reading \cite{GORV}, Chapter~2.1, as a guideline.

We assume that the number $N$ of sites is given by $N=2^K$ for some
large number $K\in\N$. The step to arbitrary $N$ is not difficult;
cf. Remark~\ref{r_arbitrary_N}, below. We decompose the spin system
into blocks, each containing two spins. The coarse-graining operator
$P\dvtx  X_{N,m} \to X_{\fraca{N}{2},m}$ assigns to each block the mean spin
of the block. More precisely, $P$ is given by
%
%e7 ###
%e7 #&#
\begin{equation}
\label{e_two_pair_coarse_grained_operator}
P(x):=  \bigl(\tfrac{1}{2} (x_1+x_2),\tfrac{1}{2} (x_3+x_4), \ldots,
\tfrac{1}{2} (x_{N-1}+x_N)  \bigr).
\end{equation}
Due to the coarse-graining operator $P$, we can decompose the canonical
ensemble $\mu_{N,m}$ into
%
%e8 ###
%e8 #&#
\begin{equation}\label{e_disintegration_example}
\mu_{N,m}(dx)= \mu(dx|y) \bar\mu(dy),
\end{equation}
where $\bar\mu:= P_{\#} \mu_{N,m}$ denotes the push forward of the
Gibbs measure $\mu$ under $P$ and $\mu(dx|y)$ is the conditional
measure of $x$ given $Px=y$. The last equation has to be understood in
a weak sense; that is, for any test function~$\xi$
\[
\int\xi\,d \mu_{N,m} = \int_Y  \biggl(\int_{ \{ Px=y  \}
} \xi \mu(dx|y)  \biggr) \bar\mu(dy).
\]
Now, we are able to state the first ingredient of the proof of
Theorem~\ref{p_main_result}.\vspace*{-1pt}

%pr2.1 #&#
\begin{proposition}[(Hierarchic criterion for the LSI)]\label{p_hierarchic}
Assume that the~sin\-gle-site potential $\psi$ is perturbed strictly
convex in the sense of \eqref{e_perturbation_convex}. If the~mar\-ginal
$\bar\mu$ satisfies the LSI with constant $\varrho_1>0$ uniformly in
the system size $N$ and the mean spin $m$, then the canonical ensemble
$\mu_{N,m}$ also satisfies the LSI with constant $\varrho_2>0$
uniformly in the system size $N$ and the mean spin $m$.\looseness=-1\vspace*{-1pt}
\end{proposition}

The proof of this statement is given in Section~\ref
{s_auxiliary_results}. Due to the last proposition it suffices to
deduce the LSI for the marginal $\bar\mu$. Hence, let us have a
closer look at the structure of $\bar\mu$. We will characterize the
Hamiltonian of the marginal $\bar\mu$ with the help of the
renormalization operator $\mathcal R$, which is introduced as follows.\vspace*{-1pt}
%
%de2.2 #&#
\begin{definition}
Let $\psi\dvtx  \R\to\R$ be a single-site potential. Then the
renormalized single-site potential $\mathcal R \psi\dvtx  \R\to\R$ is
defined by
%
%e9 ###
%e9 #&#
\begin{equation}
\label{e_defintion_renormalization_operator}
\mathcal{R} \psi(y):= - \log\int\exp \bigl( - \psi(x+y) - \psi
(-x +y)  \bigr) \,dx.\vadjust{\goodbreak}
\end{equation}
\end{definition}

%
%re2.3 #&#
\begin{remark}\label{r_interpretation_renormalization}
The renormalized single-site potential $\mathcal R \psi$ can be
interpreted in the following way: A change of variables (cf. \cite{Evans}, Section 3.3.3) and the invariance of the Hausdorff measure
under translation yield the identity
%
%e10 ###
\begin{eqnarray*}
\exp ( - \mathcal{R} \psi(y)  ) & =& \int\exp \bigl(-
\psi(x+y) - \psi(-x +y)  \bigr)\,dx \\
&=& \frac{1}{\sqrt{2}} \int\exp \bigl(- \psi(x_1) - \psi(x_2)
 \bigr) \mathcal{H}_{\lfloor{ \{x_1+x_2=2y \}} }^1 (dx).% \label{e_partition_function_vs_renormalization}
\end{eqnarray*}
Therefore, the renormalized single-site potential $\mathcal R \psi$
describes the free energy of two independent spins $X_1$ and $X_2$
[identically distributed as\break $Z^{-1} \exp(- \psi)$] conditioned on a
fixed mean value $\frac{1}{2} ( X_1+X_2  )= y$.
\end{remark}

%
%le2.4 #&#
\begin{lemma}[(Invariance under renormalization)]\label
{p_invariance_of_renormalization}
Assume that the single-site potential $\psi$ is perturbed strictly
convex in the sense of \eqref{e_perturbation_convex}. Then the
renormalized Hamiltonian $\mathcal{R} \psi$ is also perturbed
strictly convex in the sense of \eqref{e_perturbation_convex}. % More
%precisely, there is a splitting $\mathcal{R} \psi(x)= \overline{
% \[
% \overline{\psi}_c'' \gtrsim1  \mbox{and}  |\overline{
% \]
\end{lemma}

Direct calculation using the coarea formula (cf. \cite{Evans}, Section 3.4.2) reveals the following structure of the marginal $\bar
\mu$.
%
%le2.5 #&#
\begin{lemma}
The marginal $\bar\mu$ is given by
\[
\bar\mu(dy) = \frac{1}{Z}  \exp \Biggl(- \sum_{i=1}^{\fraca{N}{2}}
\mathcal{R} \psi(y_i)  \Biggr) \mathcal{H}_{\lfloor X_{\fraca
{N}{2},m} }^{\fraca{N}{2}-1}(dy).
\]
\end{lemma}

It follows from the last two lemmas that the marginal $\bar\mu$ has
the same structure as the canonical ensemble $\mu_{N,m}$. The
single-site potential of $\bar\mu$ is given by the renormalized
single-site potential $\mathcal R \psi$. Hence, one can iterate the
coarse-graining of pairs. The next statement shows that after finitely
many iterations the renormalized single-site potential $\mathcal{R}^M
\psi$ becomes uniformly strictly convex. Therefore, the Bakry--\'Emery
criterion (cf. Theorem~\ref{p_criterion_bakry_emery}) yields that the
corresponding marginal satisfies the LSI with constant $\tilde\varrho
>0$, uniformly in the system size $N$ and the mean spin $m$. Then, an
iterated application of the hierarchic criterion of the LSI
(cf. Proposition~\ref{p_hierarchic}) yields Theorem~\ref
{p_main_result} in the case $N = 2^K$.

%th2.6 #&#
\begin{theorem}[(Convexification by renormalization)]\label{p_convexification}
Let $\psi$ be a perturbed strictly convex single-site potential in the
sense of \eqref{e_perturbation_convex}. Then there is an integer $M_0$
such that for all $M \geq M_0$ the $M$-times renormalized single-site
potential $\mathcal R^M \psi$ is uniformly strictly convex
independently of the system size $N$ and the mean spin $m$.
\end{theorem}

We conclude this section by giving some remarks and pointing out the
central tools needed for the proof of the auxiliary results. The next
remark shows how Theorem~\ref{p_main_result} is verified in the case
of an arbitrary number $N$ of sites.\looseness=-1\vadjust{\goodbreak}

%re2.7 #&#
\begin{remark} \label{r_arbitrary_N}
Note that an arbitrary number of sites $N$ can be written as
\[
N = \tilde K 2^K + R
\]
for some number $\tilde K$, a large but fixed number $K$ and a bounded
number $R < 2^K$. Hence, one can decompose the spin system into $\tilde
K$ blocks of $2^K$ spins and one block of $R$ spins. The big blocks of
$2^K$ spins are coarse-grained by pairs, whereas the small block of $R$
spins is not coarse-grained at all. After iterating this procedure
sufficiently often, the renormalized single-site potentials of the big
blocks are uniformly strictly convex. On the remaining block of
$R$ spins, the corresponding single-site potentials are unchanged.
Because $\psi$ is a bounded perturbation of a strictly convex
function, it follows from a combination of the Bakry--\'Emery criterion
(cf. Theorem~\ref{p_criterion_bakry_emery}) and the Holley--Stroock
criterion (cf. Theorem~\ref{p_criterion_holley_stroock}) that the
marginal of the whole system satisfies the LSI with constant
\[
\varrho\gtrsim\exp \bigl(-R  ( \sup\delta\psi- \inf\delta
\psi )  \bigr) ,
\]
which is independent on $N$ and $m$. Therefore, an iterated application
of the hierarchic criterion of the LSI (cf. Proposition~\ref
{p_hierarchic}) yields Theorem~\ref{p_main_result}.
\end{remark}

%re2.8 #&#
\begin{remark}[(Inhomogeneous single-site potentials)]
It is a natural question whether this approach can be applied to the
case of inhomogeneous single-site potentials. In this case, the
single-site potentials are allowed to depend on the sites; that is, the
Hamiltonian has the form $H= \sum_{i=1}^N \psi_i$ where each $\psi
_i$ is a perturbed strictly-convex potential. In principle, we believe
that our approach can be adapted to this situation even if not in a
straightforward way. The reason is that only one step of the proof of
Theorem~\ref{p_main_result} has to be adapted: It is the
convexification of the single-site potentials by iterated
renormalization (see Theorem~\ref{p_convexification}).
\end{remark}

Let us make a comment on the proof of Theorem~\ref{p_convexification},
which is stated in Section~\ref{s_convexification}. Starting point for
the proof is the observation that the $M$-times renormalized
single-site potential $\mathcal R^M \psi$ corresponds to the
coarse-grained Hamiltonian related to coarse-graining with block size
$2^M$; cf. \cite{GORV}.
%
%le2.9 #&#
\begin{lemma}\label{p_renormalization_equivalent_to_coarse_graining}
For $K \in\mathbb{N}$ let the coarse-grained Hamiltonian $\bar H_K$
be defined by
%
%e11 ###
%e10 #&#
\begin{equation}\label{e_d_bar_H}
\bar H_K (m) = - \frac{1}{K} \log\int\exp(- H (x))  \mathcal
{H}_{\lfloor X_{K,m}}^{K -1} (dx).
\end{equation}
Let $M \in\mathbb{N}$. Then there is a constant $0 < C(2^M) < \infty
$ depending only on $2^M$ such that
\begin{eqnarray*}
% \label{e_renormalization_equivalent_to_coarse_graining}
\mathcal{R}^M \psi = 2^M \bar H_{2^M} + C (2^M).
\end{eqnarray*}
\end{lemma}

Because the last statement is verified by a straightforward
application of the area and coarea formula, we omit the proof.
In Lemma~\ref{p_renormalization_equivalent_to_coarse_graining} one
could easily determine the exact value of the constant
$C(2^M)$. However, the exact value is not important because we are only
interested in the convexity of
$\mathcal{R}^M \psi$. In \cite{GORV}, the convexification of $\bar
H_K$ was deduced from a local
Cram\'er theorem; cf. \cite{GORV}, Proposition~31. For the proof of
Theorem~\ref{p_convexification} we
follow the same strategy generalizing the argument to perturbed
strictly convex single-site potentials $\psi$.

Now, we make some comments on the proof of Proposition~\ref
{p_hierarchic} and Lem\-ma~\ref{p_invariance_of_renormalization}, which
are stated in Section~\ref{s_auxiliary_results}. One of the limiting
factors in the proof of Theorem~\ref{p_GORV} is the application of a
classical covariance estimate; cf.~\cite{GORV}, Lemma 22. In our
framework this estimate can be formulated as:
%
%le2.10 #&#
\begin{lemma}\label{p_cov_clas}
Assume that the single-site potential $\psi$ is perturbed strictly
convex in the sense of \eqref{e_perturbation_convex}. Let $\nu$ be a
probability measure on $\R$ given by
\[
\nu(dx) = \frac{1}{Z} \exp ( - \psi(x)  )\,dx.
\]
Then for any function $f\geq0$ and $g$
\[
| \cov_{\nu} (f,g) | \lesssim \sup_x  | g'(x)  |
\biggl( \int f \,d \nu \biggr)^{\h}   \biggl( \int\frac{| f' |^2}{f}\,d\nu
 \biggr)^{\h} .
\]
\end{lemma}

In \cite{GORV}, the last estimate was applied to the function $g= \psi
'$. Note that the function $  | g'(x)  | =  | \psi
''(x)  |$ is only bounded in the case of a perturbed quadratic
single-site potential $\psi$. The main new ingredient for the proof of
the hierarchic criterion for the LSI (cf. Proposition~\ref
{p_hierarchic}) and the invariance principle (cf. Lemma~\ref
{p_invariance_of_renormalization}) is an asymmetric Brascamp--Lieb
inequality, which does not exhibit this restriction.
%
%le2.11 #&#
\begin{lemma}\label{p_asymmetric_BL}
Assume that the single-site potential $\psi$ is perturbed strictly
convex in the sense of \eqref{e_perturbation_convex}. Let $\nu$ be a
probability measure on $\R$ given by
\[
\nu(dx) = \frac{1}{Z} \exp ( - \psi(x)  )\,dx.
\]
Then for any function $f$ and $g$
\[
| \cov_{\nu} (f,g) | \leq\exp{ ( 3 \osc\delta\psi )}
 \sup_x  \biggl| \frac{g'(x)}{\psi_c''(x)}  \biggr|  \int| f' |\,d\nu,
\]
where $\osc\delta\psi:= \sup_x \delta\psi(x) - \inf_x \delta
\psi(x)$.
\end{lemma}

We call the last inequality asymmetric because, compared to the
original Brascamp--Lieb inequality \cite{BL}, the space $L^2\times
L^2$ is replaced by $L^1\times L^\infty$, and the factor $(\psi
_c'')^{-\h}$ is not evenly distributed. It is an interesting question\vadjust{\goodbreak}
if an analog statement also holds for higher dimensions. The proof of
Lemma~\ref{p_asymmetric_BL} is based on a kernel representation of the
covariance. All steps are elementary.% We close this section with
%stating the proof of the asymmetric Brascamp-Lieb inequality.
%
\begin{pf*}{Proof of Lemma~\ref{p_asymmetric_BL}}
Let us consider a Gibbs measure $\mu$ associated to the Hamiltonian
$H\dvtx  \R\to\R$. More precisely, $\mu$ is given by
\[
\mu(dx) := \frac{1}{Z} \exp ( - H(x)  )\,dx .
\]
We start by deriving the following integral representation of the
covariance of $\mu$:
%
%e12 ###
%e11 #&#
\begin{equation}\label{L.1}
\operatorname{cov}_\mu(f,g) =
\int\int f'(x)K_{\mu}(x,y) g'(y)\,dx\,dy,
\end{equation}
where the nonnegative kernel $K_{\mu}(x,y)$ is given by
\[
K_{\mu}(x,y) :=\left\{
\begin{array}{l@{ \quad }l}
M_{\mu}(x)(1-M_{\mu})(y)&\mbox{for } y\ge x\\
(1-M_{\mu})(x)M_{\mu}(y)&\mbox{for } y\le x
\end{array}
 \right\},
\]
and $M_{\mu}(x) := \mu((-\infty,x))$ so that $(1-M_{\mu})(x)=\mu
((x,\infty))$.
Indeed, we start by noting that
%
%e13 ###
%e12 #&#
\begin{equation}\label{L.3}
 \qquad    \operatorname{cov}_{\mu}(f,g) = \int\int\bigl(f(z)-f(x)\bigr)\mu(x)\,dx
\int\bigl(g(z)-g(y)\bigr)\mu(y)\,dy \,\mu(z)\,dz,
\end{equation}
where we do not distinguish between the measure $\mu(dx)$ and its
Lebesgue density $\mu(x)$ in our notation.
Using $M_{\mu}'(x)=\mu(x)$, we can use integration by parts
to rewrite each factor in terms of the derivative
\begin{eqnarray}
 &&\int\bigl(f(z)-f(x)\bigr)\mu(x)\,dx \nonumber\\
&& \qquad =\int_{-\infty} ^z\bigl(f(z)-f(x)\bigr) M_{\mu}'(x)\,dx
-\int_{z}^{\infty}\bigl(f(z)-f(x)\bigr)(1-M_{\mu})'(x)\,dx\nonumber\\
&& \qquad =\int_{-\infty}^z f'(x) M_{\mu}(x)\,dx
-\int_{z}^{\infty} f'(x) (1-M_{\mu})(x)\,dx\nonumber\\
&& \qquad =\int f'(x) \bigl(I(x<z)M_{\mu}(x)-I(x>z)(1-M_{\mu})(x) \bigr)\,
dx,\nonumber
\end{eqnarray}
where $I(x<z)$ assumes the value $1$ if $x<z$ and zero otherwise.
Inserting this and the corresponding identity for $g(y)$ into (\ref
{L.3}), we obtain
%
%e13 #&#
\begin{eqnarray}
\operatorname{cov}_{\mu}(f,g)
&=&\int\! \int f'(x) \bigl(I(x<z)M_{\mu}(x) - I(x>z)(1-M_{\mu})(x)
\bigr)\,dx\nonumber\\
& & \hphantom{\int\! \int}{}\times\int g'(y) \bigl(I(y<z)M_{\mu}(y)-I(y>z)(1-M_{\mu})(y)
\bigr)\,dy \mu(z)\,dz\\
& =&\int\int f'(x) K_{\mu}(x,y) g'(y)\,dx\,dy
\nonumber
\end{eqnarray}
with kernel $K_{\mu}(x,y)$ as desired, given by
\begin{eqnarray}
&&K_{\mu}(x,y)\nonumber\\
&& \qquad = M_{\mu}(x)M_{\mu}(y){\int} I(x<z)I(y<z)\mu(z)\,
dz\nonumber\\
&& \qquad  \quad {} -M_{\mu}(x)(1-M_{\mu})(y){\int} I(x<z)I(y>z)\mu(z)\,
dz\nonumber\\
&& \qquad  \quad {} -(1-M_{\mu})(x)M_{\mu}(y){\int} I(x>z)I(y<z)\mu(z)\,
dz\nonumber\\
&& \qquad  \quad {} +(1-M_{\mu})(x)(1-M_{\mu})(y){\int} I(x>z)I(y>z)\mu
(z)\,dz\nonumber\\
&& \qquad = M_{\mu}(x)M_{\mu}(y)(1-M_{\mu})(\max\{x,y\})\nonumber\\
&& \qquad  \quad {} -M_{\mu}(x)(1-M_{\mu})(y)I(y>x)\bigl(M_{\mu}(y)-M_{\mu}(x)\bigr)\nonumber
\\
&& \qquad  \quad {} -(1-M_{\mu})(x)M_{\mu}(y)I(y<x)\bigl(M_{\mu}(x)-M_{\mu}(y)\bigr)\nonumber
\\
&& \qquad  \quad {} +(1-M_{\mu})(x)(1-M_{\mu})(y)M_{\mu}(\min\{x,y\})\nonumber\\
&& \qquad =I(y>x)
 \bigl(M_{\mu}(x)M_{\mu}(y)(1-M_{\mu})(y) \nonumber\\
&& \qquad    \hphantom{=I(y>x)
 \bigl(}{}- M_{\mu}(x)(1-M_{\mu})(y)\bigl(M_{\mu}(y)-M_{\mu}(x)\bigr)\nonumber\\
&& \hspace*{24.5pt}\qquad    \hphantom{=I(y>x)
 \bigl(}{}+(1-M_{\mu})(x)(1-M_{\mu})(y)M_{\mu}(x) \bigr)\nonumber\\
&& \qquad  \quad {} +I(y\le x)
 \bigl(M_{\mu}(x)M_{\mu}(y)(1-M_{\mu})(x) \nonumber\\
&& \qquad  \quad\hphantom{{} +I(y\le x)
 \bigl(} {}- (1-M_{\mu})(x)M_{\mu}(y)\bigl(M_{\mu}(x)-M_{\mu}(y)\bigr)\nonumber\\
&& \hspace*{25pt}\qquad  \quad\hphantom{{} +I(y\le x)
 \bigl(} {}+(1-M_{\mu})(x)(1-M_{\mu})(y)M_{\mu}(y) \bigr)\nonumber\\
&& \qquad =I(y>x)M_{\mu}(x)(1-M_{\mu})(y)+I(y\le x)(1-M_{\mu})(x)M_{\mu
}(y).\nonumber
\end{eqnarray}

We now establish the following identity for the above kernel:
%
%e14 ###
%e14 #&#
\begin{equation}\label{L.2}
\int K_{\mu}(x,y) H''(y)\,dy = \mu(x).
\end{equation}
Indeed, we have by integrations by part
\begin{eqnarray}
&&\int K_{\mu}(x,y) H''(y)\,dy\nonumber\\
&& \qquad =(1-M_{\mu})(x)\int_{-\infty}^{x}M_{\mu}(y) H''(y)\,dy
+ M_{\mu}(x) \int_{x}^{\infty}(1-M_{\mu})(y) H''(y)\,dy\nonumber\\
&& \qquad =(1-M_{\mu})(x)\biggl (M_{\mu}(x) H'(x)
-\int_{-\infty}^x M_{\mu}'(y)H'(y)\,dy \biggr)
\nonumber\\
&& \qquad  \quad {} +M_{\mu}(x) \biggl(-(1-M_{\mu})(x) H'(x)
+\int_x^\infty M_{\mu}'(y) H'(y)\,dy \biggr)\nonumber\\
&& \qquad =-(1-M_{\mu})(x)\int_{-\infty} ^x\exp(-H(y)) H'(y)\,dy\nonumber\\
&& \qquad  \quad {} + M_{\mu}(x) \int_{x}^{\infty}\exp(-H(y)) H'(y)\,dy\nonumber\\
&& \qquad =(1-M_{\mu})(x)\mu(x)+M_{\mu}(x)\mu(x) = \mu(x).\nonumber
\end{eqnarray}
Let us now consider the Gibbs measures $\nu(dx)$ and $\nu_c(dx)$,
given by
\[
\nu(dx) = \frac{1}{Z} \exp \bigl( - \psi_c (x) - \delta\psi(x)
 \bigr)\,dx   \quad \mbox{and} \quad \nu_c(dx) = \frac{1}{Z} \exp
 ( - \psi_c (x)  )\,dx.
\]
By the integral representation \eqref{L.1} of the covariance we have
the estimate
\[
 | \operatorname{cov}_\nu(f,g)  |  \leq
\int\int | f'(x) | K_{\nu}(x,y)  | g'(y)  | \,
dx\,dy.
\]
By a straight-forward calculation, we can estimate
\begin{eqnarray*}
M_{\nu}(x) & =& \frac{\int_{- \infty}^x \exp(- \psi_c (x) - \delta
\psi(x))\,dx} {\int\exp(- \psi_c (x) - \delta\psi(x))\,dx } \\
& \leq&\exp( \osc\delta\psi)   \frac{\int_{- \infty}^x \exp(-
\psi_c (x) )\,dx} {\int\exp(- \psi_c (x) )\,dx } \\
& =& \exp( \osc\delta\psi)   M_{\nu_c} (x).
\end{eqnarray*}
Together with a similar estimate for $  (1 - M_{\nu}(y)
),$ this yields the kernel estimate
\[
K_{\nu}(x,y) \leq\exp(2 \osc\delta\psi)   K_{\nu_c}(x,y).
\]
Applying this to the covariance estimate from above yields
\begin{eqnarray*}
 | \cov_\nu(f,g)  |  \leq  \exp(2 \osc\delta\psi)
\int\int | f'(x) | K_{\nu_c}(x,y)  | g'(y)  |
\,dx\,dy.
\end{eqnarray*}
Using the identity \eqref{L.2} for $\mu= \nu_c$, we may easily conclude
\begin{eqnarray*}
 | \cov_\nu(f,g)  |   &\leq&  \exp(2 \osc\delta\psi)
  \sup_y \frac{ | g'(y) | }{\psi_c'' (y)}   \int |
f'(x) | \int K_{\nu_c}(x,y) \psi_c'' (y) \,dy\, dx \\
& =&   \exp(2 \osc\delta\psi)   \sup_y \frac{ | g'(y) | }{\psi
_c'' (y)}   \int | f'(x) | \nu_c(dx) \\
& \leq&   \exp(3 \osc\delta\psi)   \sup_y \frac{ | g'(y) |
}{\psi_c'' (y)}   \int | f'(x) | \nu(dx) .
\end{eqnarray*}
\upqed
\end{pf*}

For the entertainment of the reader, let us argue how the identity
(\ref{L.2}) also yields the traditional Brascamp--Lieb inequality in
the case $H''>0$. Indeed, by the symmetry of the kernel $K_\mu(x,y)$,
identity \eqref{L.2} yields, for all $x$ and $y$,
%
%e15 ###
%e15 #&#
\begin{equation}\label{e_kernel_identity}
\int K_\mu(x,y)H''(y)\,dy = \mu(x)\quad\mbox{and} \quad\int
K_\mu(x,y)H''(x)\,dx = \mu(y).
\end{equation}
The integral representation of the covariance (\ref{L.1}) yields
\begin{eqnarray*}
\var_{\mu} (f) & =&   \int\int f'(x)K_{\mu}(x,y) f'(y)\,dx\,dy \\
& =&   \int\int f'(x)  \biggl( \frac{K_{\mu}(x,y)   H''(y)}{H''(x)}
 \biggr)^\h   f'(y)  \biggl( \frac{K_{\mu}(x,y)   H''(x)}{H''(y)}
 \biggr)^\h \,dx\,dy.
\end{eqnarray*}
Then a combination of H\"{o}lder's inequality and the identity \eqref
{e_kernel_identity} for the kernel $K_\mu(x,y)$ yields the
Brascamp--Lieb inequality,
%
%e16 ###
%e16 #&#
\begin{eqnarray} \label{e_BL}
\var_{\mu} (f) &  \leq &  \biggl( \int\int \frac
{|f'(x)|^2}{H''(x)}   K_{\mu}(x,y) H''(y)\,dy\,dx  \biggr)^\h \nonumber\\
&&{}   \times\biggl ( \int\int \frac{|f'(y)|^2}{H''(y)}
K_{\mu}(x,y) H''(x)\,dx\,dy  \biggr)^\h \nonumber
\\[-8pt]
\\[-8pt]
&  = &   \biggl( \int \frac{|f'(x)|^2}{H''(x)} \mu(x)\,dx  \biggr)^\h
   \biggl( \int \frac{|f'(y)|^2}{H''(y)} \mu(y)\,dy  \biggr)^\h\nonumber\\
&  = &  \int \frac{|f'(x)|^2}{H''(x)} \mu(x)\,dx.
\nonumber
\end{eqnarray}

%s2.2 ###
%s2.2 #&#
\subsection{Proof of auxiliary results}\label{s_auxiliary_results}
In this section we outline the proof of Proposition~\ref{p_hierarchic}
and Lemma~\ref{p_invariance_of_renormalization}. We start with
Proposition~\ref{p_hierarchic}, which is the hierarchic criterion for
the LSI. Unfortunately, we cannot directly apply the two-scale
criterion of \cite{GORV}, Theorem 3. The reason is that the number
%
%e17 ###
%e17 #&#
\begin{equation}
\label{e_interaction_between_scales}\quad
\kappa:= \max \biggl\{ \frac{\skp{ \Hess H (x) u}{v}}{|u|  |v|},
u \in\im(2P^tP),  v \in\im(\id_X - 2 P^tP)  \biggr\},
\end{equation}
which measures the interaction between the microscopic and macroscopic
scales, can be infinite for a perturbed strictly convex single-site
potential $\psi$. However, we follow the proof of \cite{GORV}, Theorem
3, with only one major difference: Instead of applying the classical
covariance estimate (cf. Lemma~\ref{p_cov_clas}), we apply the
asymmetric Brascamp--Lieb inequality; cf. Lemma~\ref{p_asymmetric_BL}.
Let us assume for the rest of this section that the single-site
potential $\psi$ is perturbed strictly convex in the sense of \eqref
{e_perturbation_convex}.

For convenience, we set $X:= X_{N,m}$ and $Y:=X_{\fraca{N}{2},m}$. We
choose on $X$ and $Y$ the standard Euclidean structure given by
\[
\skp{x}{y} = \sum_{i=1}^N x_i y_i.
\]
The coarse-graining operator $P\dvtx X \to Y$ given by \eqref
{e_two_pair_coarse_grained_operator} satisfies the identity
\[
% \label{e_identity_coarse_graining_operator}
2 PP^t= id_Y,
\]
where $P^t\dvtx Y \to X$ is the adjoint operator of $P$. Note that our $P^t$
differs from the $P^t$ of \cite{GORV}, because the Euclidean structure
on Y differs from the Euclidean structure used in \cite{GORV} by a
factor. The last identity yields that $2 P^t P$ is the orthogonal
projection of $X$ to $\im P^t$. Hence, one can decompose $X$ into the
orthogonal sum of \emph{microscopic fluctuations} and \emph
{macroscopic variables} according to
\[
X  = \ker P \oplus\im P^t
\]
and
\[
x  =  ( \id_X - 2 P^t P  ) x + 2 P^t P x.
\]
We apply this decomposition to the gradient $\nabla f$ of a smooth
function $f$ on $X$. The gradient $\nabla f$ is decomposed into a
macroscopic gradient and a fluctuation gradient satisfying
%
%e19 ###
%e18 ###
%e18 #&#
\begin{eqnarray}
\label{e_decomposition_of_gradients}
\nabla f (x) & =&  ( \id_X - 2 P^t P  ) \nabla f(x) + 2 P^t
P \nabla f(x) \quad\mbox{and} \nonumber
\\[-8pt]
\\[-8pt]
 | \nabla f (x)  |^2 & =&  | ( \id_X - 2 P^t P
 ) \nabla f(x)  |^2 +  | 2 P^t P \nabla f(x)  |^2.
\nonumber
\end{eqnarray}
Note that $\ker P$ is the tangent space of the fiber $ \{Px=y
 \}$. Hence the gradient of $f$ on $ \{Px=y  \}$ is
given by $ ( \id_X - 2 P^t P  ) \nabla f(x)$. The first
main ingredient of the proof of Proposition~\ref{p_hierarchic} is the
following statement.
%
%le2.12 #&#
\begin{lemma} \label{p_micro_LSI}
The conditional measure $\mu(dx|y)$ given by \eqref
{e_disintegration_example} satisfies the LSI with constant $\varrho>0$
uniformly in the system size $N$, the macroscopic profile $y$ and the
mean spin $m$. More precisely, for any nonnegative function~$f$
\begin{eqnarray*}
&&\int f \log f \mu(dx|y) - \int f \mu(dx|y) \log \biggl( \int f \mu
(dx|y)  \biggr) \\
&& \qquad \leq\frac{1}{2 \varrho} \int\frac{| ( \id_X - 2 P^t P
) \nabla f|^2}{f} \mu(dx|y).
\end{eqnarray*}
\end{lemma}
\begin{pf*}{Proof of Lemma~\ref{p_micro_LSI}}
Observe that the conditional measure $\mu(dx| y)$ has a product
structure: We decompose $ \{ Px=y  \}$ into a product of
Euclidean spaces. Namely for
\[
X_{2,y_i}:=  \{ (x_{2i-1}, x_{2i} ),   x_{2i-1} + x_{2i} = 2 y_i
 \}, \qquad i \in\biggl \{1, \ldots, \frac{N}{2}  \biggr\},
\]
we have
\[
 \{ Px=y  \} = X_{2, y_1} \times\cdots\times X_{2,
y_{\fraca{N}{2}}}.
\]
It follows from the coarea formula (cf. \cite{Evans}, Section 3.4.2) that
\begin{eqnarray*}
&& \int_{ \{ Px=y  \}} f(x) \mu(dx|y) \\
&& \qquad   = \int f(x) \bigotimes_{i=1 }^{\fraca{N}{2}} \frac{1}{Z} \exp
 \bigl(- \psi(x_{2i-1}) - \psi(x_{2i})  \bigr)  \mathcal
{H}_{\lfloor X_{2,y_i}}^{1} (dx_{2i-1}, d x_{2i}).
\end{eqnarray*}
Hence $\mu(dx|y)$ is the product measure
%
%e20 ###
%e19 #&#
\begin{equation} \label{e_product_mu}
\mu(dx|y) = \bigotimes_{i=1}^{\fraca{N}{2}} \mu_{2,y_i}(dx_{2i -1},
dx_{2 i} ),
\end{equation}
where we make use of the notation introduced in \eqref
{e_definition_of_mu_N_m}. Because the single-site potential $\psi$ is
perturbed strictly convex in the sense of \eqref
{e_perturbation_convex}, a combination of the Bakry--\'Emery criterion
(cf. Theorem~\ref{p_criterion_bakry_emery}) and the Holley--Stroock
criterion (cf. Theorem~\ref{p_criterion_holley_stroock}) yield that
the measure $\mu_{2,m}(dx_{1}, dx_{2})$ satisfies the LSI with
constant $\varrho>0$ uniformly in $m$. Then the tensorization
principle (cf. Theorem~\ref{p_tensorization_principle}) implies the
desired statement.
\end{pf*}

For convenience, let us introduce the following notation: Let $f$ be an
arbitrary function. Then its conditional expectation $\bar f$ is
defined by
\[
% \label{e_notation_coarse_grained_function}
\bar f (y) := \int f(x) \mu(dx| y).
\]
The second main ingredient of the proof of Proposition~\ref
{p_hierarchic} is the following proposition, which is the analog
statement of \cite{GORV}, Proposition~20.
%
%pr2.13 #&#
\begin{proposition}\label{p_estimate_of_macroscopic_fisher_information}
Assume that the marginal $\bar\mu(dy)$ given by \eqref
{e_disintegration_example} satisfies the LSI uniformly in the system
size $N$ and the mean spin $m$. Then for any nonnegative function $f$,
\begin{eqnarray*}
\frac{|\nabla\bar f (y)|^2}{\bar f(y) } \lesssim \int\frac
{|\nabla f (x)|^2}{f(x)}  \mu(dx| y),
\end{eqnarray*}
uniformly in the macroscopic profile $y$ and the system size $N$.
\end{proposition}

Before we verify Proposition~\ref
{p_estimate_of_macroscopic_fisher_information}, let us show how it can
be used in the proof of Proposition~\ref{p_hierarchic}.
\begin{pf*}{Proof of Proposition~\ref{p_hierarchic}}
Using Lemma~\ref{p_micro_LSI} and Proposition~\ref
{p_estimate_of_macroscopic_fisher_information} from above, the argument
is exactly the same as in the proof of \cite{GORV},  Theorem~3:

Let $\phi$ denote the function $\phi(x):= x \log x$. The additive
property of the entropy implies
\begin{eqnarray*}
\int\phi(f)\,d\mu_{N,m} - \phi \biggl( \int f d \mu_{N,m}  \biggr)
&=& \int \biggl[ \int\phi (f(x)  ) \mu(dx|y) - \phi
( \bar f(y)  )  \biggr] \bar\mu(dy) \\
&&{}  +  \biggl[ \int\phi (\bar f (y)  ) \bar \mu(dy)
- \phi \biggl( \int\bar f(y) \bar\mu(dy)  \biggr)  \biggr] .
\end{eqnarray*}
An application of Lemma~\ref{p_micro_LSI} yields the estimate
\begin{eqnarray*}
&&\int \biggl[ \int\phi (f(x)  ) \mu(dx|y) - \phi (
\bar f(y)  )  \biggr] \bar\mu(dy) \\
&& \qquad \leq\frac{1}{2 \varrho} \int\int \frac{| ( \id_X - 2 P^tP
 ) \nabla f (x) |^2}{f(x)} \mu(dx| y) \bar\mu(dy).
\end{eqnarray*}
By assumption the marginal $\bar\mu$ satisfies the LSI with constant
$\lambda>0$. Together with Proposition~\ref
{p_estimate_of_macroscopic_fisher_information} this yields the estimate
\begin{eqnarray*}
\int\phi (\bar f (y)  ) \bar \mu(dy) - \phi\biggl ( \int
\bar f(y) \bar\mu(dy)  \biggr)
& \leq&\frac{1}{2 \lambda}
\int \frac{|\nabla\bar f(y)|^2}{\bar f(y)}  \bar\mu(dy) \\
& \lesssim&\int\int\frac{|\nabla f (x)|^2}{f(x)}  \mu(dx| y) \bar
\mu(dy) .
\end{eqnarray*}
A combination of the last three formulas and the observations \eqref
{e_disintegration_example} and \eqref{e_decomposition_of_gradients} yield
\begin{eqnarray*}
&& \int\phi(f)\,d\mu_{N,m} - \phi\biggl ( \int f d \mu_{N,m}  \biggr)
\\
&& \qquad   \lesssim \int\frac{| ( \id_X - 2 P^tP  ) \nabla
f (x) |^2}{f(x)} \mu_{N,m} (dx) + \int\frac{|\nabla f (x)|^2}{f(x)}
\mu_{N,m} (dx) \\
&& \qquad   \lesssim\int\frac{|\nabla f (x)|^2}{f(x)} \mu_{N,m} (dx),
\end{eqnarray*}
uniformly in the system size $N$ and the mean spin $m$.
\end{pf*}

Because the hierarchic criterion for the LSI is an important ingredient
in the proof of the main result, we outline the proof of
Proposition~\ref{p_estimate_of_macroscopic_fisher_information} in full
detail. We follow the proof of \cite{GORV}, Proposition~20, which is
based on two lemmas. We directly take over the first lemma (cf. \cite{GORV}, Lemma~21), which in our notation becomes:
%
%le2.14 #&#
\begin{lemma} \label{p_decomposition_of_mean_of_microscopic_gradient}
For any function $f$ on $X$ and any $y \in Y$, it holds
\[
% \label{e_decomposition_of_mean_of_microscopic_gradient}
\int P \nabla f(x) \mu(dx|y) = \frac{1}{2} \nabla\bar f (y) + P \cov
_{\mu(dx|y)} (f, \nabla H).
\]
\end{lemma}

%
%re2.15 #&#
\begin{remark}
The notational difference compared to \cite{GORV}, Lemma~21, is based
on our choice of the Euclidean structure on $Y=X_{\fraca{N}{2},m}$.
Compared to the notation in Lemma~21 of \cite{GORV}, we have
\[
\nabla_Y \bar f(y)= \frac{N}{2} \nabla\bar f(y).
\]
Hence we omit the proof, which is a straightforward calculation.
\end{remark}

The more interesting ingredient of the proof of \cite{GORV}, Proposition~20, is the estimate (see \cite{GORV}, (42), (43))
\[
\bigl|2P \cov_{\mu(dx|y)} (f, \nabla H)\bigr|^2 \leq \frac{\sqrt{2} \kappa
^2}{\varrho^2}  \bar f (y)   \int\frac{| (\id_X - 2 P^t P) \nabla
f(x) |^2}{f(x)} \mu(dx|y) .  %
\]
In \cite{GORV}, the last estimate is deduced by direct calculation
from the standard covariance estimate given by Lemma~\ref{p_cov_clas}.
In contrast to \cite{GORV} we cannot use this estimate because the
constant $\kappa$ given by \eqref{e_interaction_between_scales} may
be infinite for a perturbed strictly convex single-site potential $\psi
$. We avoid this problem by applying the more robust asymmetric
Brascamp--Lieb inequality given by Lemma~\ref{p_asymmetric_BL}. Our
substitute for the last estimate is:
%
%le2.16 #&#
\begin{lemma}\label{p_adapted_two_scale_key_estimate}
For any nonnegative function $f$
\begin{eqnarray}
\bigl|2P \cov_{\mu(dx|y)} (f, \nabla H)\bigr|^2  \lesssim \bar f (y)   \int
\frac{| \nabla f(x) |^2}{f(x)} \mu(dx|y), \nonumber%
\end{eqnarray}
uniformly in the system size $N$, the macroscopic profile $y$ and the
mean spin~$m$.
\end{lemma}

We postpone the proof of Lemma~\ref{p_adapted_two_scale_key_estimate}
and show how it is used in the proof of Proposition~\ref
{p_estimate_of_macroscopic_fisher_information} (cf. proof of \cite
{GORV}, Proposition~20).
\begin{pf*}{Proof of Proposition~\ref{p_estimate_of_macroscopic_fisher_information}}
Note that because for any $a,b \in\R$,
\[
\tfrac{1}{2} (a + b)^2 \leq a^2 + b^2,
\]
it follows form the definition \eqref
{e_two_pair_coarse_grained_operator} of $P$ that for any $x$,
%
%e21 ###
%e20 #&#
\begin{equation} \label{e_operator_norm_P}
|Px|^2 \leq\tfrac{1}{2}  |x|^2 .
\end{equation}
%
% Note that by \eqref{e_identity_coarse_graining_operator}, we have for
%any $x \in X$
% \begin{equation}
% \label{e_identity_norm_of_macroscopic_profile}
% |2 P^t Px|^2 = 2 |Px|^2.
% \end{equation}
By successively using Lemma~\ref
{p_decomposition_of_mean_of_microscopic_gradient} and Jensen's
inequality (with the convex function $(a,b) \mapsto |b|^2 / a$), we have
%
%e22 ###
\begin{eqnarray}
\frac{|\nabla\bar f (y)|^2}{\bar f(y) } & =& \frac{4}{\bar f(y)}
 \biggl| P \int \nabla f(x) \mu(dx|y) - P \cov_{\mu(dx|y)} (f,
\nabla H)  \biggr|^2 \nonumber\\
& \lesssim&\frac{1}{\bar f(y)}  \biggl| \int P \nabla f(x) \mu
(dx|y) \biggr|^2 + \frac{1}{\bar f(y)}  \bigl| P \cov_{\mu(dx|y)}
(f, \nabla H)  \bigr|^2 \nonumber \\
& \lesssim& \int \frac{| P \nabla f(x)|^2}{f(x)} \mu(dx|y) + \frac
{1}{\bar f(y)}  \bigl| 2 P \cov_{\mu(dx|y)} (f, \nabla H)  \bigr|^2.\nonumber
 % \label{e_estimation_macroscopic_fisher_proof_1}
\end{eqnarray}
On the first term on the r.h.s. we apply the estimate \eqref
{e_operator_norm_P}. On the second term we apply Lemma~\ref
{p_adapted_two_scale_key_estimate}, which yields the desired estimate.
\end{pf*}

Now, we prove Lemma~\ref{p_adapted_two_scale_key_estimate}, which also
represents one of the main differences compared to the two-scale
approach of \cite{GORV}. The main ingredients are the product
structure (\ref{e_product_mu}) of $\mu(dx|y)$ and the asymmetric
Brascamp--Lieb inequality; cf. Lemma~\ref{p_asymmetric_BL}.
\begin{pf*}{Proof of Lemma~\ref{p_adapted_two_scale_key_estimate}}
We have to estimate the covariance
%
%e23 ###
%e21 #&#
\begin{equation} \label{e_decomp_cov}
\bigl|2P \cov_{\mu(dx|y)} (f, \nabla H)\bigr|^2 = \sum_{j=1}^{\fraca{N}{2}}
\bigl|\cov_{\mu(dx|y)}  ( f, (2P \nabla H)_j  )\bigr |^2 .
\end{equation}
Therefore, let us consider for $j \in \{1,\ldots,\frac{N}{2}
 \}$ the term $\cov_{\mu(dx|y)}  (f, (2P \nabla H
)_j)$. Note that the function
\[
\label{e_macroscopic_hamiltonian_component}
 (2 P \nabla H (x) )_j = \psi'(x_{2j-1}) + \psi'(x_{2j}) %
%=: g(x_{2j-1}, x_{2_j} ),
\]
only depends of the variables $x_{2j-1}$ and $x_{2j}$. Hence, the
product structure \eqref{e_product_mu} of $\mu(dx|y)$ yields the identity
%
%e24 ###
%e22 #&#
\begin{eqnarray} \label{e_cov_prod}\quad
&& \cov_{\mu(dx | y)} (f,2  (P \nabla H  )_j) \nonumber
\\[-8pt]
\\[-8pt]
&& \qquad    = \int\cov_{\mu_{2,y_j}(dx_{2j-1}, dx_{2_j})} (f,  (2
P \nabla H  )_j) \bigotimes_{i=1, i \neq j}^{\fraca{N}{2}} \mu
_{2, y_i}(dx_{2i-1}, dx_{2_i}).
\nonumber
\end{eqnarray}
As we will show below, we obtain, by using the asymmetric
Brascamp--Lieb inequality of Lemma~\ref{p_asymmetric_BL} and the
Csisz\'ar--Kullback--Pinsker inequality, the estimate
%
%e25 ###
%e23 #&#
\begin{eqnarray} \label{e_estimation_covariance_component}
&&  \bigl| \cov_{\mu_{2,y_j}(dx_{2j-1}, dx_{2_j})} (f,  (2 P
\nabla H  )_j)  \bigr|\nonumber\\
&& \qquad  \lesssim  \biggl( \int f(x) \mu_{2,
y_j}(dx_{2j-1}, dx_{2_j})  \biggr)^{\h}\nonumber
\\[-8pt]
\\[-8pt]
&&  \quad \qquad {}  \times  \biggl( \int \frac{| (\fracc{d}{dx_{2j-1}}) f(x) |^2 +
| (\fracc{d}{dx_{2j}}) f(x) |^2 }{f (x)}\nonumber\\
&&\hspace*{90.5pt}\hphantom{\quad \qquad {}  \times  \biggl( \int}{}\times \mu_{2, y_j}(dx_{2j-1},
dx_{2_j})  \biggr)^{\h}
\nonumber
\end{eqnarray}
uniformly in $j$ and $y_j$. Therefore, a combination of identity \eqref
{e_cov_prod}, the last estimate and H\"{o}lder's inequality yield
\begin{eqnarray}
&& \bigl| \cov_{\mu(dx | y)} (f, (2P \nabla H)_j )\bigr|^2 \nonumber\\
&& \qquad   \lesssim\int f (x) \mu(dx|y) \int \frac{|(\fracc
{d}{dx_{2j-1}}) f(x) |^2 + | (\fracc{d}{dx_{2j}}) f(x) |^2 }{f (x)} \mu
(dx|y), \nonumber% \label{e_estimation_covariance_component_2}
\end{eqnarray}
which implies the desired estimate by the identity \eqref
{e_decomp_cov}.\vadjust{\goodbreak}

It is only left to deduce estimate \eqref
{e_estimation_covariance_component}. We assume w.l.o.g. $j=1$. Recall
the splitting $\psi=\psi_c + \delta\psi$ given by \eqref
{e_perturbation_convex}. We use the bound on $|\delta\psi'|$ to estimate
%
%e26 ###
%e24 #&#
\begin{eqnarray}\label{p_split_covariance}
 &&\bigl| \cov_{\mu_{2,y_1}(dx_{1}, dx_{2})} (f, (2P \nabla H)_1 )
 \bigr| \nonumber\\
 &&  \qquad \lesssim  \bigl| \cov_{\mu_{2. y_1}(dx_{1}, dx_{2} )}
 \bigl( f, \psi_c'(x_{1}) + \psi_c'(x_{2})  \bigr)  \bigr|\\
&& \qquad  \quad {}   + \int \biggl| f -\int f \mu_{2, y_1}(dx_{1}, dx_{2} )  \biggr|
\mu_{2, y_1}(dx_{1}, dx_{2} ).
\nonumber
\end{eqnarray}
Now, we consider the first term on the r.h.s. of the last estimate. For
$y_1 \in\R$ let the one-dimensional probability measure $\nu(dz|
y_1)$ be defined by the density
%
%e27 ###
%e25 #&#
\begin{equation}\label{e_d_nu}
\nu(dz| y_1) := \frac{1}{Z}  \exp \bigl(- \psi(-z + y_1) - \psi(
z + y_1)  \bigr) \,dz.
\end{equation}
A reparametrization of the one-dimensional Hausdorff measure implies
%
%e28 ###
%e26 #&#
\begin{equation}\label{e_repara}
\int\xi(x_{1} , x_{2}) \mu_{2,y_1} (dx_{1}, d x_{2}) = \int \xi(-z
+ y_1 , z +y_1 ) \nu(dz| y_1)
\end{equation}
for any measurable function $\xi$. We may assume w.l.o.g. that $f(x)=
f(x_{1}, x_{2})$ just depends on the variables $x_{1}$ and $x_{2}$.
Hence for
\[
\tilde f (z,y_1) := f (-z +y_1, z + y_1) \quad\mbox{and} \quad
\tilde g (z, y_1) := \psi_c'(-z + y_1) + \psi_c'(z + y_1),
\]
the last identity yields
\[
\cov_{\mu_{2, y_1}(dx_{1}, dx_{2})}  \bigl( f, \psi_c'(x_{1}) + \psi
_c'(x_{2})  \bigr) = \cov_{\nu(dz | y_1)} ( \tilde f, \tilde g ).
\]
Because
\[
 \biggl| \frac{(\fracc{d}{dz})\tilde g(z,y_1)}{\psi_c''(-z + y_1) + \psi
_c''(z + y_1)}  \biggr| =  \biggl| \frac{ -\psi_c'' (-z+y_1) + \psi
_c'' (z+y_1)}{\psi_c'' (-z+y_1) + \psi_c'' (z+y_1)}  \biggr| \leq2,
\]
an application of the asymmetric Brascamp--Lieb inequality
(cf. Lemma~\ref{p_asymmetric_BL}) yields
\begin{eqnarray*}
 \bigl| \cov_{\nu(dz | y_1)} ( \tilde f, \tilde g ) \bigr | &
\lesssim& \int\biggl| \frac{d}{dz} \tilde f \biggr|   \nu(dz | y_1) \\
&\lesssim& \biggl (\int \tilde f  \nu(dz | y_1)  \biggr)^\h  \biggl(
\int\frac{| (\fracc{d}{dz}) \tilde f |^2}{ \tilde f}  \nu(dz | y_1)
 \biggr)^\h.
\end{eqnarray*}
From the last inequality and from \eqref{e_repara} follows the estimate
%
%e29 ###
%e27 #&#
\begin{eqnarray} \label{e_using_aBL}
&&  \bigl| \cov_{\mu_{2,y_1}(dx_{1}, dx_{2} )}  \bigl( f, \psi
_c'(x_{1}) + \psi_c'(x_{2})  \bigr)  \bigr| \nonumber
\\
&& \qquad   \lesssim \biggl(\int f  \mu_{2, y_1}(dx_{1} , dx_{2} )
\biggr)^\h \\
&& \qquad  \quad {}\times\biggl ( \int\frac{| (\fracc{d}{dx_{1}}) f |^2 + | (\fracc
{d}{dx_{2}}) f |^2 }{ f}  \mu_{2,y_1}(dx_{1} , dx_{2})  \biggr)^\h.
\nonumber
\end{eqnarray}
We turn to the second term on the r.h.s. of \eqref
{p_split_covariance}. For convenience, let us write $\bar f (y_1) : =
\int f \mu_{2,y_1}(dx_{1}, dx_{2} )$. An application of the well-known
Csisz\'ar--Kullback--Pinsker inequality (cf. \cite{Csi,Kul}) yields
\begin{eqnarray*}
&&\int | f - \bar f(y_1)  | \mu_{2, y_1}(dx_{1}, dx_{2})\\
 && \qquad  =
\bar f(y_1)  \int \biggl| \frac{f}{\bar f (y_1)} - 1  \biggr| \mu
_{2,y_1}(dx_{1}, dx_{2})\\
&& \qquad  \lesssim\bar f (y_1)  \biggl ( \int\frac{f}{\bar f (y_1)} \log
\frac{f}{\bar f (y_1)}  \mu_{2,y_1}(dx_{1}, dx_{2})  \biggr)^{\h}.
\end{eqnarray*}
An application of the LSI for the measure $\mu_{2,y_1}(dx_{1}, dx_{2}
)$ implies (cf. proof of Lemma~\ref{p_micro_LSI})
\begin{eqnarray*}
&& \int \biggl| f-\int f \mu_{2, y_1}(dx_{1}, dx_{2} )  \biggr|  \mu
_{2, y_1}(dx_{1}, dx_{2}) \\
&& \qquad   \lesssim\biggl ( \int f \mu_{2,y_1}(dx_{1}, dx_{2} )
\biggr)^{\h} \\
&& \qquad  \quad {}\times \biggl( \int\frac{| (\fracc{d}{dx_{1}}) f |^2 + | (\fracc
{d}{dx_{2}}) f |^2}{f} \mu_{2, y_1}(dx_{1}, dx_{2})  \biggr)^{\h}.
\end{eqnarray*}
A combination of \eqref{p_split_covariance}, \eqref{e_using_aBL}, and
the last inequality yield the estimate \eqref
{e_estimation_covariance_component}.
\end{pf*}

We turn to the proof of Lemma~\ref{p_invariance_of_renormalization}.
Again, the main ingredient of the proof is the asymmetric
Brascamp--Lieb inequality.
\begin{pf*}{Proof of Lemma~\ref{p_invariance_of_renormalization}}
We define
\[
% \label{e_bar_psi_c}
\overline{\psi}_c (m) := - \frac{1}{2} \log\int\exp \bigl(- \psi
_c (-z+ m) - \psi_c  ( z + m  )  \bigr)\,dz
\]
and
\begin{eqnarray*}
\overline{\delta\psi} (m) &: =& - \frac{1}{2} \log\int\exp
\bigl(- \psi(-z+ m) - \psi ( z + m  )  \bigr)\,dz \\
&&{}   + \frac{1}{2} \log\int\exp \bigl(- \psi_c (-z+ m) - \psi
_c  ( z + m  )  \bigr)\,dz. % \label{e_bar_delta_psi}
\end{eqnarray*}
Now, we show that the splitting $\mathcal{R} \psi= \overline{\psi
}_c + \overline{\delta\psi}$ satisfies the conditions given
by \eqref{e_perturbation_convex}. Using the strict convexity of $\psi
_c$ it follows by a standard argument based on the Brascamp--Lieb
inequality (cf. \cite{BL} and \eqref{e_BL}) that the first condition
is preserved, that is,
\[
\overline{\psi}_c'' \gtrsim1.
\]

We turn to the perturbation $\overline{\delta\psi}$. Analogously to
the measure $\nu(dz|m)$ given by \eqref{e_d_nu}, we introduce the
measure $\nu_c(dz|m)$ via the density
\[
\nu_c (dz) := \frac{1}{Z}   \exp \bigl( - \psi_c (-z+ m) - \psi_c
 ( z + m  )  \bigr)\,dz.
\]
It follows that
\[
\overline{\delta\psi} (m) = - \frac{1}{2} \log\int\exp \bigl( -
\delta\psi(-z+ m) - \delta\psi ( z + m  )  \bigr) \nu
_c (dz).
\]
Direct calculation using the bound $|\delta\psi| \lesssim1$ yields
\begin{eqnarray*}
|\overline{\delta\psi} (m)|  \lesssim1.
\end{eqnarray*}

We turn to the first derivative of $\overline{\delta\psi}$. A direct
calculation based on the definition of $\overline{\delta\psi}$ yields
\begin{eqnarray*}
2 \overline{\delta\psi}'(m) & =& \int \bigl(\psi' (-z+ m) + \psi'
 ( z + m  )  \bigr) \nu(dz) \\
&&{}  - \int\bigl (\psi_c' (-z+ m) + \psi_c'  ( z + m
)  \bigr) \nu_c (dz).
\end{eqnarray*}
For $s \in [ 0, 1  ]$ we define the measure $\nu^s (dz)$
by the probability density
\begin{eqnarray*}
\frac{1}{Z}   \exp \bigl( - \psi_c (-z+ m) - \psi_c  ( z + m
 ) - s \delta\psi(-z+ m) - s \delta\psi ( z + m  )
 \bigr)\,dz.
\end{eqnarray*}
Note that $\nu^s$ interpolates between $\nu^0= \nu_c$ and $\nu^1 =
\nu$. By the mean-value theorem there is $s \in [ 0,1  ]$
such that
\begin{eqnarray*}
&&\hspace*{-3pt} 2 \overline{\delta\psi}' (m) \\
&&\hspace*{-3pt} \qquad  = \frac{d}{ds} \int \bigl(\psi_c' (-z+ m) + \psi_c'  ( z + m
 ) + s \delta\psi' (-z+ m) + s \delta\psi'  ( z + m
 )  \bigr) \nu^s (dz) \\
&&\hspace*{-3pt} \qquad  = \int \bigl( \delta\psi' (-z+ m) + \delta\psi'  ( z + m
 )  \bigr) \nu^{s} (dz) \\
&&\hspace*{-3pt} \qquad  \quad {}  + \cov_{\nu^{s}} \bigl ( \psi_c' (-z+ m) + \psi_c'  (
z + m  ) ,  \delta\psi(-z+ m) + \delta\psi ( z + m
 )  \bigr) \\
&&\hspace*{-3pt} \qquad  \quad {}  + \cov_{\nu^{s}} \bigl ( s \delta\psi' (-z+ m) + s \delta
\psi'  ( z + m  ) ,   \delta\psi(-z+ m) + \delta\psi
 ( z + m  )  \bigr).
\end{eqnarray*}
The first term on the r.h.s. is controlled by the assumption $|\delta
\psi'| \lesssim1$. We turn to the estimation of the first covariance
term. An application of the asymmetric Brascamp--Lieb inequality of
Lemma~\ref{p_asymmetric_BL} and $|\delta\psi| + |\delta\psi'|
\lesssim1$ yields the estimate
\begin{eqnarray*}
&&  \bigl| \cov_{\nu^{s}}  \bigl( \psi_c' (-z+ m) + \psi_c'  (
z + m  ) ,  \delta\psi(-z+ m) + \delta\psi ( z + m
 )  \bigr)  \bigr| \\
&& \qquad  \lesssim \sup_z  \biggl| \frac{\psi_c'' (-z+ m) - \psi_c''
( z + m  )}{\psi_c'' (-z+ m) + \psi_c''  ( z + m  )}
 \biggr|  \\
 && \qquad  \quad {}\times \int | - \delta\psi' (-z+ m) + \delta\psi'
( z + m  )  | \nu^{s} (dz) \\
&& \qquad  \lesssim1.
\end{eqnarray*}
The second covariance term can be estimated by using $|\delta\psi| +
|\delta\psi'| \lesssim1$. Summing up, we have deduced the desired
estimate $ | \overline{\delta\psi}' | \lesssim1$.
\end{pf*}

%s3 ###
%s3 #&#
\section{Convexification by iterated renormalization}\label{s_convexification}
In this section we prove Theorem~\ref{p_convexification} that states
the convexification of a perturbed strictly convex single-site
potential $\psi$ by iterated renormalization. The proof relies on a
local Cram\'er theorem and some auxiliary results. The proof of
Theorem~\ref{p_convexification} is given in Section~\ref
{s_s_convexification}. The proofs of the auxiliary results are given in
Section~\ref{s_local_cramer}.

%s3.1 ###
%s3.1 #&#
\subsection{\texorpdfstring{Proof of Theorem~\protect\ref{p_convexification}}
{Proof of Theorem 2.6}}\label{s_s_convexification}

Let us consider the coarse-grained Hamiltonian $\bar H_K$ given
by \eqref{e_d_bar_H}. In view of Lemma~\ref
{p_renormalization_equivalent_to_coarse_graining}, it suffices to show
the strict convexity of $\bar H_K$ for large $K\gg1$. The strategy is
the same as in \cite{GORV}, Proposition 31. Let $\varphi$ denote the
Cram\'er transform of $\psi$, namely
\[
\varphi(m) := \sup_{\sigma\in\R}\biggl  ( \sigma m - \log\int\exp
\bigl(\sigma x - \psi(x)\bigr)\,dx  \biggr).
\]
Because $\varphi$ is the Legendre transform of the strictly convex function
%
%e30 ###
%e28 #&#
\begin{equation}
\label{e_definition_varphi*}
\varphi^* (\sigma) = \log\int\exp\bigl(\sigma x - \psi(x)\bigr)\,dx,
\end{equation}
there exists for any $m \in\R$, a unique $\sigma=\sigma(m)$, such that
%
%e31 ###
%e29 #&#
\begin{equation}
\label{e_identity_legendre_transform}
\varphi(m) = \sigma m - \varphi^*(\sigma).
\end{equation}
From basic properties of the Legendre transform, it follows that
$\sigma$ is determined by the equation
%
%e32 ###
%e30 #&#
\begin{equation}
\label{e_dsigma_phi*}
\frac{d}{d \sigma} \varphi^* (\sigma) = \frac{\int x \exp(\sigma
x - \psi(x))\,dx}{\int\exp(\sigma x - \psi(x))\,dx} = m.
\end{equation}
The starting point of the proof of the convexification of the
coarse-grained Hamiltonian $\bar H_K(m)$ is the explicit representation
%
%e33 ###
%e31 #&#
\begin{equation}
\label{e_cramer_idea}
\tilde g_{K,m}(0)= \exp\bigl (K \varphi(m) - K \bar H_K (m)  \bigr).
\end{equation}
Here, $\tilde g_{K,m}$ denotes the Lebesgue density of the distribution
of the random variable
\[
\frac{1}{\sqrt{K}} \sum_{i=1}^K  ( X_i - m  ),
\]
where $X_i$ are $K$ real-valued independent random variables
identically distributed according to
%
%e34 ###
%e32 #&#
\begin{equation}
\label{e_definition_mu_sigma}
\mu^{\sigma}(dx) := \exp \bigl( - \varphi^* (\sigma) + \sigma x -
\psi(x)  \bigr)\,dx.
\end{equation}
We note that in view of \eqref{e_dsigma_phi*} the mean of $X_i$ is
$m$. As in \cite{GORV}, (125), the Cram\'er representation \eqref
{e_cramer_idea} follows from direct substitution and the coarea
formula. As we will see in the proof of Lemma~\ref
{p_crucial_convexification}, the Cram\'er transform $\varphi$ is
strictly convex. The main idea of the proof is to transfer the
convexity from $\varphi$ to $\bar H_K$ using representation \eqref
{e_cramer_idea} and a local central limit type theorem for the density
$\tilde g_{K,m}$, which is formulated in the next statement. %Note that
%the statement of Proposition~\ref{P} is not restricted to perturbed
%quadratic Hamiltonians as was the local Cram\'er theorem used in

%pr3.1 #&#
\begin{proposition}\label{P}
Let $\psi(x)$ be a smooth function that is increasing sufficiently
fast as $|x|\uparrow\infty$ for all subsequent integrals to exist.
Note that the probability measure $\mu^\sigma$ defined by \eqref
{e_definition_mu_sigma} depends on the field strength $\sigma$. We
introduce its mean $m$ and variance $s^2$
%
%e35 ###
%e33 #&#
\begin{equation}\label{1}
m :=  \int x \mu^\sigma(dx)\quad\mbox{and}\quad
s^2 :=  \int(x-m)^2 \mu^\sigma(dx) .
\end{equation}

We assume that uniformly in the field strength $\sigma$,
the probability measure $\mu^\sigma$ has its standard deviation
$s$ as unique length scale in the sense that
%
%e37 ###
%e36 ###
%e34 #&#
%e35 #&#
\begin{eqnarray}\label{4}
\int|x-m|^k \mu^\sigma(dx)&\lesssim&s^k\qquad\mbox{for }
k=1,\ldots,5,\\\label{7}
 \biggl| \int\exp(ix\xi)\mu^\sigma(dx)  \biggr|&\lesssim&|s\xi
|^{-1} \qquad\mbox{for all } \xi\in\R .
%&s^2\langle(\frac{df}{dx})^2\rangle\quad\mbox{for all}\;f(x).\label{4}
\end{eqnarray}

Consider $K$ independent random variables $X_1,\ldots,X_K$ identically
distributed according to $\mu^\sigma$. Let $g_{K,\sigma}$ denote
the Lebesgue density of the distribution of the normalized sum
$\frac{1}{\sqrt{K}}\sum_{i=1}^K\frac{X_i-m}{s}$.

Then $g_{K,\sigma}(0)$ converges for $K\uparrow\infty$
to the corresponding value for the normalized Gaussian.
This convergence is uniform in $m$, of order $\frac{1}{\sqrt{K}}$, and
$C^2$ in $\sigma$:
%
%e40 ###
%e39 ###
%e38 ###
%e36 #&#
%e37 #&#
%e38 #&#
\begin{eqnarray}\label{9}
\biggl|g_{K,\sigma}(0)-\frac{1}{\sqrt{2\pi}}\biggr|&\lesssim&\frac{1}{\sqrt
{K}},\\\label{11}
\biggl|\frac{1}{s}\frac{d}{d\sigma}g_{K, \sigma}(0)\biggr|&\lesssim&\frac
{1}{\sqrt{K}},\\\label{12}
\biggl|\biggl(\frac{1}{s}\frac{d}{d\sigma}\biggr)^2g_{K, \sigma}(0)\biggr|&\lesssim&
\frac{1}{\sqrt{K}}.
\end{eqnarray}
%
%(Note that $\frac{1}{s}\frac{d}{d\sigma}$ is the nondimensional
%version of the derivative w.\ r.\ t.\ the parameter).
\end{proposition}

Let us comment a bit on this result: Quantitative versions of
the central limit theorem like (\ref{9}) are abundant in the literature;
see, for instance,~\cite{Feller}, Chapter XVI, \cite
{KipLand}, Appendix~2, \cite{GPV}, Section~3, and \cite{LPY}, page~752
and Section~5. In his work on the spectral gap,
Caputo appeals even to a finer estimate that
makes the first terms in an error expansion in $K^{- \h}$ explicit~\cite{Cap},
Theorem~2.1. The coefficients of the higher order terms
are expressed in terms of moments of $\mu^\sigma$.
However, following \cite{GORV}, Proposition~31, for our two-scale
argument we need
\textit{pointwise} control of the Lebesgue density $g_{K , \sigma}$
[in form of $g_{K, \sigma}(0)$] and, in addition, control of derivatives
of $g_{K, \sigma}$ w.r.t. the field parameter $\sigma$; cf.
(\ref{11}), (\ref{12}). Note that the derivative $\frac{d}{d\sigma}$
has units of length (because $\sigma$, which multiplies $x$ in the
Hamiltonian [cf. (\ref{e_definition_mu_sigma})] has units of inverse
length) so that
$\frac{1}{s}\frac{d}{d\sigma}$ is the properly nondimensionalized
derivative. Pointwise control means that control of the moments
[cf. (\ref{4})] is not sufficient. One also needs to know that
$\mu^\sigma$
has no fine structure on scales much smaller than~$s$.
This property is ensured the upper bound (\ref{7}).

As opposed to \cite{GORV}, Proposition~31, the Hamiltonian
$\psi$ we want Proposition~\ref{P} to apply
is not a perturbation of the quadratic $\frac{1}{2}x^2$, but of
a general, strictly convex potential $\psi$. As a consequence, the
variance $s^2$
can be a strongly varying function of the field strength $\sigma$.
Nevertheless, Lemma~\ref{p_assumptions_lct} from below shows that
every element
$\mu^\sigma$ in the family of measures
is characterized by the single length scale~$s$, uniformly
in $\sigma$ in the sense of (\ref{4}) and (\ref{7}). For the
verification of (\ref{4}) in Lemma~\ref{p_assumptions_lct}, one could
take over the argument of \cite{Cap}, Lemma~2.2, that relies on a
result by Bobkov \cite{Bob} stating that the SG constant $\varrho$ of
the measure $\mu^\sigma$ can be estimated by its variance, that
is, $\varrho\gtrsim\frac{1}{s^2}$. However, we provide a
self-contained argument for the verification of (\ref{4}) and (\ref
{7}) in Lemma~\ref{p_assumptions_lct} just using basic calculus of one
variable. The merit of Proposition~\ref{P} consists in providing
a version of the central limit theorem that is $C^2$ in
the field strength $\sigma$ even if the variance $s^2$ varies
strongly with $\sigma$.
%
%le3.2 #&#
\begin{lemma}\label{p_assumptions_lct}
Assume that the single-site potential $\psi$ is perturbed strictly
convex in the sense of \eqref{e_perturbation_convex}. Then $s \lesssim
1$ uniformly in $m$, and conditions \eqref{4} and \eqref{7} of
Proposition~\ref{P} are satisfied.
\end{lemma}

Using Proposition~\ref{P}, Lemma~\ref{p_assumptions_lct}, and the
Cram\'er representation \eqref{e_cramer_idea} we could easily deduce a
local Cram\'er theorem (cf. \cite{GORV}, Proposition 31) for general
perturbed strictly convex potentials $\psi$. However, because we are
just interested in the convexification of $\bar H_K$, we just consider
the convergence of the second derivatives of $\varphi$ and $\bar H_K$.
%
%le3.3 #&#
\begin{lemma}\label{p_crucial_convexification}
Assume that the single-site potential $\psi$ is perturbed strictly
convex in the sense of \eqref{e_perturbation_convex}. Then for all $m
\in\R$ it holds
\[
 \biggl| \frac{d^2}{dm^2}  \varphi(m) -   \frac{d^2}{dm^2}  \bar
H_K (m)  \biggr| \lesssim\frac{1}{K s^2},
\]
where $s^2$ is defined as in Proposition~\ref{P}.
\end{lemma}

\begin{pf*}{Proof of Theorem~\ref{p_convexification}}
Because of Lemma~\ref{p_renormalization_equivalent_to_coarse_graining}
it suffices to show that there exists $\delta> 0 $ and $K_0 \in\N$
such that for all $K \geq K_0$ and $m \in\R$
\[
\frac{d^2}{dm^2} \bar H_K(m) \geq\delta.
\]
We start with some formulas on the derivatives of $\varphi$.
Differentiation of identity \eqref{e_identity_legendre_transform} yields
\begin{eqnarray*}
\frac{d}{dm} \varphi& =& \frac{d}{dm} \sigma  m + \sigma - \frac
{d}{d\sigma} \varphi^*  \frac{d}{dm} \sigma\\
& \overset{\mathclap{\eqref{e_dsigma_phi*}}}{=}& \frac{d}{dm} \sigma
 m + \sigma - m  \frac{d}{dm} \sigma\\
& =& \sigma.
\end{eqnarray*}
A direct calculation reveals that [see \eqref{e_d_dsigma_m} below]
\begin{eqnarray*}
\frac{d}{d \sigma} m  = s^2,
\end{eqnarray*}
where $s^2$ is defined as in Proposition~\ref{P}. Hence, a second
differentiation of $\varphi$ yields the identity
%
%e41 ###
%e39 #&#
\begin{equation} \label{e_formula_dd_dmm_phi}
\frac{d^2}{dm^2} \varphi= \frac{d}{dm} \sigma = \biggl (\frac{d}{d
\sigma} m  \biggr)^{-1} = \frac{1}{s^2}.
\end{equation}
By Lemma~\ref{p_crucial_convexification} we thus have
\begin{eqnarray*}
\frac{d^2}{dm^2} \bar H_K & =& \frac{d^2}{dm^2} \varphi+ \frac
{d^2}{dm^2}  (\bar H_K - \varphi  ) \\
& \geq&\frac{1}{s^2} - \frac{C}{K} \frac{1}{s^2} \\
& \geq&\frac{1}{2}   \frac{1}{s^2},
\end{eqnarray*}
if $K \geq K_0$ for some large $K_0$. The statement follows from the
uniform bound $s \lesssim1$ provided by Lemma~\ref{p_assumptions_lct}.
\end{pf*}

%s3.2 ###
%s3.2 #&#
\subsection{Proof of the local Cram\'er theorem and of the auxiliary
results}\label{s_local_cramer}

In this section we prove the auxiliary statements of the last
subsection. Before turning to the proof of Proposition~\ref{P} we
sketch the strategy. For convenience we introduce the notation
%
%e42 ###
%e40 #&#
\begin{equation}\label{e_notation_integral}
\langle f\rangle := \int f (x) \mu^\sigma(dx)= \int f(x)\exp
\bigl(-\varphi^*(\sigma)+\sigma x-\psi(x)\bigr)\,dx.
\end{equation}
The definition of $g_{K, \sigma}$ (cf. Proposition~\ref{P}) suggests
to introduce the shifted and rescaled variable
%
%e43 ###
%e41 #&#
\begin{equation}\label{5}
\hat x := \frac{x-m}{s}.
\end{equation}
We note that by (\ref{1}) the first and second moment in $\hat x$
are normalized
%
%e44 ###
%e42 #&#
\begin{equation}\label{2}
\langle\hat x\rangle = 0,\qquad\langle\hat x^2\rangle = 1
\end{equation}
and that (\ref{4}) turns into
%Like for all central limit theorems, one key ingredient for
%Proposition~\ref{P} are higher moment bounds that
%are uniform in $\sigma$. Those follow easily from the spectral gap
%estimate (\ref{4}):
%
%Under the assumptions of Proposition~\ref{P} we have
%
%e45 ###
%e43 #&#
\begin{equation}\label{3}
\sum_{k=1}^5\langle|\hat x|^k\rangle \lesssim 1.
\end{equation}

Proposition~\ref{P} is a version of the central limit theorem that,
like most others, is best proved with help of the Fourier transform.
Indeed, since the random variables $\hat X_1:=\frac{X_1-m}{s},\ldots,
\hat X_K:=\frac{X_K-m}{s}$ in the statement of
Proposition~\ref{P} are independent and identically distributed, the
distribution of their sum is the $K$-fold convolution of
the distribution of $\hat X_1$. Therefore, the Fourier transform
of the distribution of the $\sum_{n=1}^K\hat X_n$ is the $K$th power
of the Fourier transform of the distribution of $\hat X$.
The latter is given by
\[
\langle\exp(i\hat x\hat\xi)\rangle,
\]
where $\hat\xi$ denotes the variable dual to $\hat x$. Hence,\vspace*{1pt} the Fourier
transform of the distribution of the normalized sum
$\frac{1}{\sqrt{K}}\sum_{n=1}^K\hat X_K$\vspace*{2pt} is given by
$\langle\exp(i\hat x\frac{1}{\sqrt{K}}\hat\xi)\rangle^K$. Applying
the inverse Fourier transform, we obtain the representation
%
%e46 ###
%e44 #&#
\begin{equation}\label{13bis}
2\pi g_{K, \sigma}(0) =
\int\biggl\langle\exp\biggl(i\hat x\frac{1}{\sqrt{K}}\hat\xi\biggr)\biggr\rangle^K\,d\hat
\xi.
\end{equation}

In order to make use of formula (\ref{13bis}), we need estimates
on $\langle\exp(i\hat x\hat\xi)\rangle$.
Because of
%
%e47 ###
%e45 #&#
\begin{eqnarray}\label{19}
\frac{d^k}{d\hat\xi^k}\langle\exp(i\hat x\hat\xi)\rangle
 =
i^k\langle\hat x^k\exp(i\hat x\hat\xi)\rangle,
\end{eqnarray}
the moment bounds (\ref{3}) translate into control of
$\langle\exp(i\hat x\hat\xi)\rangle$ for $|\hat\xi|\ll1$.
Together with the normalization (\ref{2}), we obtain, in particular,
\[
\bigl|\langle\exp(i\hat x\hat\xi)\rangle-\bigl(1-\tfrac{1}{2}\hat\xi^2\bigr)\bigr|
\lesssim
 |\hat\xi|^3.
\]
We will use the latter in the following form: There exists
a complex-valued function $h(\hat\xi)$ such that for $|\hat\xi|\ll1$,
%
%e48 ###
%e46 #&#
\begin{equation}\label{C}
\langle\exp(i\hat x\hat\xi)\rangle = \exp(-h(\hat\xi))\qquad
\mbox{with }
\bigl|h(\hat\xi)-\tfrac{1}{2}\hat\xi^2\bigr| \lesssim |\hat\xi|^3.
\end{equation}
This estimate, showing that the Fourier transform of the normalized
probability $\langle\cdot\rangle$ is close for $|\hat\xi|\ll1$
to the Fourier transform
of the normalized Gaussian, is at the core of most proofs of the
central limit
theorem.

Estimate \eqref{C} provides good control over
$\langle\exp(i\hat x\hat\xi)\rangle$ for $|\hat\xi|\ll1$.
Another key ingredient is uniform decay for
$|\hat\xi|\gg1$. In our new variables, (\ref{7}) takes on the form
%
%Under the assumptions of Proposition~\ref{P} we have
%
%e49 ###
%e47 #&#
\begin{equation}\label{L2}
|\langle\exp(i\hat x\hat\xi)\rangle| \lesssim
  |\hat\xi|^{-1}.\vadjust{\goodbreak}
\end{equation}

As usual in central limit theorems, we also need control of the
characteristic function for intermediate values of $|\hat\xi|$.
This can be inferred from (\ref{3}) and (\ref{L2}) by a soft argument
(in particular, it
does not require the more intricate argument for \cite{Cap}, (2.10),
from \cite{Cap}, Lemma 2.5):
%
%le3.4 #&#
\begin{lemma}\label{soft}
Under the assumptions of Proposition~\ref{P} and for any
\mbox{$\delta>0$}, there exists $\lambda<1$ such that for all $\sigma$,
\[
|\langle\exp(i\hat x\hat\xi)\rangle| \le \lambda\qquad
\mbox{for all } |\hat\xi|\ge\delta.
\]
\end{lemma}

So far, the strategy is standard; now comes the new ingredient:
In view of formula (\ref{13bis}), in order to control $\sigma$-derivatives
of $g_{K, \sigma}(0)$, we need to control
$\frac{1}{s}\frac{d}{d\sigma}\langle\exp(i\hat x\hat\xi)\rangle$.
Relying on the identities
%
%e51 ###
%e50 ###
%e48 #&#
%e49 #&#
\begin{eqnarray}\label{14}
\frac{1}{s}\frac{d}{d\sigma}\langle f(x)\rangle&=&\langle\hat x
f(x)\rangle,
\\
\label{15}
\frac{1}{s}\frac{d}{d\sigma}\hat x&=&-1-\frac{1}{2}\langle\hat
x^3\rangle\hat x
\end{eqnarray}
that will be established in the proof of Lemma~\ref{L3} below,
we see that the estimate again follows from the moment control (\ref{3}).
Lemma~\ref{L3} is the only new element of our analysis.

%le3.5 #&#
\begin{lemma}\label{L3}
Under the assumptions of Proposition~\ref{P} we have
%
%e53 ###
%e52 ###
%e50 #&#
%e51 #&#
\begin{eqnarray}\label{21}
\biggl|\frac{1}{s}\frac{d}{d\sigma}\langle\exp(i\hat x\hat\xi)\rangle\biggr|
&\lesssim&
(1+|\hat\xi|)|\hat\xi|^3,\\\label{20}
\biggl|\biggl(\frac{1}{s}\frac{d}{d\sigma}\biggr)^2\langle\exp(i\hat x\hat\xi
)\rangle\biggr|
&\lesssim&
(1+\hat\xi^2)|\hat\xi|^3.
\end{eqnarray}
\end{lemma}

%%%%%%%%%%%%%%%%%%%%%%%%%%%%%%%%%%%%%%%%%%%%%%%%%%%%%%%%%%%%%%%%%%%%%

Before we deduce Proposition~\ref{P}, we prove Lemma~\ref{soft} and
Lemma~\ref{L3}.

\begin{pf*}{Proof of Lemma~\ref{soft}}
In view of (\ref{3}) and (\ref{L2}), it suffices to show: For
any $C<\infty$ and $\delta>0$ there exists $\lambda<1$ with
the following property: Suppose $\langle\cdot\rangle$ is a probability
measure (in $\hat x$) such that
%
%e55 ###
%e54 ###
%e52 #&#
%e53 #&#
\begin{eqnarray}\label{s1}
\langle|\hat x|\rangle&\le&C,\\\label{s2}
|\langle\exp(i\hat x\hat\xi)\rangle|&\le&\frac{C}{|\hat\xi
|}\qquad
\mbox{for all } \hat\xi.
\end{eqnarray}
Then
\[
|\langle\exp(i\hat x\hat\xi)\rangle| \le \lambda
\qquad\mbox{for all } |\hat\xi|\ge\delta.
\]
In view of (\ref{s2}), it is enough to show
\[
|\langle\exp(i\hat x\hat\xi)\rangle| \le \lambda
\qquad\mbox{for all } \delta\le|\hat\xi|\le\frac{1}{\delta}.\vadjust{\goodbreak}
\]

We give an indirect argument for this statement and
thus assume that there is
a sequence $\{\langle\cdot\rangle_{\nu}\}$ of probability measures
satisfying (\ref{s1}) and (\ref{s2}) and a sequence $\{\hat\xi_\nu\}$
of numbers in $[\delta,\frac{1}{\delta}]$ such that
%
%e56 ###
%e54 #&#
\begin{equation}\label{s5}
\liminf_{\nu\uparrow\infty}|\langle\exp(i\hat x\hat\xi_\nu
)\rangle_\nu| \ge 1.
\end{equation}
In view of (\ref{s1}), after passage to a subsequence, we may assume
that there exists a~probability measure $\langle\cdot\rangle_{\infty}$
and a number $\hat\xi_\infty>0$ such that
%
%e58 ###
%e57 ###
%e55 #&#
%e56 #&#
\begin{eqnarray}\label{s3}
\lim_{\nu\uparrow\infty}\langle f\rangle_\nu&=&\langle f\rangle
_\infty\qquad
\mbox{for all bounded and continuous } f(\hat x),\\\label{s4}
\lim_{\nu\uparrow\infty}\hat\xi_\nu&=&\hat\xi_\infty.
\end{eqnarray}
Since $|\exp(i\hat x\hat\xi_\nu)-\exp(i\hat x\hat\xi_\infty
)|\le
|\hat x||\hat\xi_\nu-\hat\xi_{\infty}|$, we obtain the following
from (\ref{s1}),
(\ref{s3}) and (\ref{s4}):
\[
\lim_{\nu\uparrow\infty}\langle\exp(i\hat x\hat\xi_\nu)\rangle
_\nu =
\langle\exp(i\hat x\hat\xi_\infty)\rangle_\infty,
\]
so that (\ref{s5}) saturates to
%
%e59 ###
%e57 #&#
\begin{equation}\label{s6}
|\langle\exp(i\hat x\hat\xi_\infty)\rangle_\infty| \ge 1.
\end{equation}

On the other hand, (\ref{s2}) is preserved under (\ref{s3})
so that we have, in particular,
%
%e60 ###
%e58 #&#
\begin{equation}\label{s7}
\lim_{|\hat\xi|\uparrow\infty}|\langle\exp(i\hat x\hat\xi
)\rangle_\infty| = 0.
\end{equation}
We claim that (\ref{s6}) and (\ref{s7}) contradict each other.
Indeed, since $\hat x\mapsto\break\exp(i\hat x\hat\xi_\infty)$ is
$S^1$-valued, it follows from (\ref{s6}) that there is a fixed
$\zeta\in S^1$ such that
\[
\exp(i\hat x\hat\xi_\infty) = \zeta\qquad
\mbox{for } \langle\cdot\rangle_\infty\mbox{-a.e. }\hat x.
\]
This implies for every $n\in\mathbb{N}$,
\[
\exp(i\hat x (n\hat\xi_\infty)) = \zeta^n\qquad
\mbox{for } \langle\cdot\rangle_\infty\mbox{-a.e. }\hat x
\]
and thus
%
%e61 ###
%e59 #&#
\begin{equation}
|\langle\exp(i\hat x (n\hat\xi_\infty))\rangle_\infty| =
|\zeta^n| = 1,
\end{equation}
which, in view of $\hat\xi_\infty\not=0$ and
thus $|n\hat\xi_\infty|\uparrow\infty$ as $n\uparrow\infty$,
contradicts (\ref{s7}).
\end{pf*}

%%%%%%%%%%%%%%%%%%%%%%%%%%%%%%%%%%%%%%%%%%%%%%%%%%%%%%%%%%%%%%%%

\begin{pf*}{Proof of Lemma~\ref{L3}}
We restrict our attention to estimate (\ref{20});
estimate (\ref{21}) is easier and can be derived by the same arguments.
We start with the identities (\ref{14}) and (\ref{15}). % Deriving (
% w.r.t.~$\sigma$ yields, making use of the abbreviation~(\ref{1}):
% %
% \begin{equation}\nonumber
% \frac{d\varphi^*}{d\sigma} = -m.
% \end{equation}
%
Deriving \eqref{e_notation_integral} w.r.t. $\sigma$ yields
%
%e62 ###
%e60 #&#
\begin{equation}\label{17}
\frac{d}{d\sigma}\langle f(x)\rangle =
\biggl\langle\biggl(x -\frac{d\varphi^*}{d\sigma} \biggr)f(x)\biggr\rangle  \overset
{\mathclap{\eqref{e_dsigma_phi*}}}{=}
\langle(x-m)f(x)\rangle.
\end{equation}
In view of definition (\ref{5}), the latter turns into (\ref{14}).\vadjust{\goodbreak}

We now turn to identity (\ref{15}) and note that, in view of
definitions (\ref{1}) and (\ref{5}), the identity (\ref{17}) yields,
in particular,
%
%e64 ###
%e63 ###
%e61 #&#
%e62 #&#
\begin{eqnarray}\label{e_d_dsigma_m}
\frac{d}{d\sigma}m&\stackrel{\mathclap{(\ref{1}),(\ref{17})}}{=}&
\langle(x-m)x\rangle \stackrel{\mathclap{(\ref{1})}}{=} \langle
(x-m)^2\rangle
 \stackrel{\mathclap{(\ref{1})}}{=} s^2,\label{36}\\
\label{e_d^2_d^2sigma_m}
\frac{d}{d\sigma}s^2 & \stackrel{\mathclap{(\ref{1}),(\ref{17})}}{=}&
\langle(x-m)(x-m)^2\rangle \stackrel{\mathclap{(\ref{5})}}{=} s^3\langle
\hat x^3\rangle,
\end{eqnarray}
which we rewrite as
\begin{eqnarray}
\frac{1}{s}\frac{d}{d\sigma}m&=&s, \nonumber\\
\frac{1}{s}\frac{d}{d\sigma}s&=&\frac{1}{2}s\langle\hat x^3\rangle
. \nonumber
\end{eqnarray}
These formulas imply, as desired,
\[
\frac{1}{s}\frac{d}{d\sigma}\hat x \stackrel{\mathclap{(\ref{5})}}{=}
\frac{1}{s}\frac{d}{d\sigma}\frac{x-m}{s} =
-1-\frac{1}{2}\langle\hat x^3\rangle\hat x.
\]

We now combine formulas \eqref{14} and (\ref{15}) to
express derivatives of $\langle f(\hat x)\rangle$. We start
with the first derivative,
%
%e65 ###
%e63 #&#
\begin{eqnarray}
\label{18}
\frac{1}{s}\frac{d}{d\sigma}\langle f(\hat x)\rangle
&\stackrel{\mathclap{(\ref{14})}}{=}&
\biggl\langle\frac{df}{d\hat x}(\hat x)\frac{1}{s}\frac{d}{d\sigma}\hat x
+f(\hat x)\hat x\biggr\rangle\nonumber
\\[-8pt]
\\[-8pt]
&\stackrel{\mathclap{(\ref{15})}}{=}&
-\biggl\langle\frac{df}{d\hat x}(\hat x)\biggr\rangle
-\frac{1}{2}\langle\hat x^3\rangle\biggl\langle\hat x \frac{df}{d\hat
x}(\hat x)\biggr\rangle
+\langle\hat x f(\hat x)\rangle.
\nonumber
\end{eqnarray}
[As a consistency check we note that
$\frac{1}{s}\frac{d}{d\sigma}\langle f(\hat x)\rangle\stackrel
{\mathclap{(\ref{18})}}{=}
-\langle(\frac{d}{d\hat x}-\hat x)f\rangle
-\frac{1}{2}\langle\hat x^3\rangle
\langle\hat x \frac{df}{d\hat x}\rangle$
vanishes if $\psi$ is quadratic since then the distribution
of $\hat x$ under $\langle\cdot\rangle$ is the normalized Gaussian
so that both $\langle(\frac{d}{d\hat x}-\hat x)f\rangle=0$
and $\langle\hat x^3\rangle=0$.]

Iterating this formula, we obtain for the second derivative,
\begin{eqnarray*}
&&\biggl(\frac{1}{s}\frac{d}{d\sigma}\biggr)^2\langle f(\hat x)\rangle\\%
&&\qquad{}\stackrel{\mathclap{(\ref{18})}}{=}
-\frac{1}{s}\frac{d}{d\sigma}\biggl\langle\frac{df}{d\hat x}(\hat
x)\biggr\rangle
-\frac{1}{2} \biggl(\frac{1}{s}\frac{d}{d\sigma}\langle\hat
x^3\rangle \biggr)
\biggl\langle\hat x \frac{df}{d\hat x}(\hat x)\biggr\rangle\\
&&\qquad\quad{}-\frac{1}{2}\langle\hat x^3\rangle
 \biggl(\frac{1}{s}\frac{d}{d\sigma}\biggl\langle\hat x \frac{df}{d\hat
x}(\hat x)\biggr\rangle \biggr)
+\frac{1}{s}\frac{d}{d\sigma}\langle\hat x f(\hat x)\rangle\\
&&\qquad{}\stackrel{\mathclap{(\ref{18})}}{=}
\biggl\langle\frac{d^2f}{d\hat x^2}\biggr\rangle
+\frac{1}{2}\langle\hat x^3\rangle
\biggl\langle\hat x \frac{d^2f}{d\hat x^2}\biggr\rangle
-\biggl\langle\hat x \frac{df}{d\hat x}\biggr\rangle\\
&&\qquad\quad{}+ \frac{1}{2} \biggl(3\langle\hat x^2\rangle+\frac{3}{2}\langle\hat x^3\rangle^2
-\langle\hat x^4\rangle \biggr)
\biggl\langle\hat x \frac{df}{d\hat x}\biggr\rangle\\
&&\qquad\quad{}+\frac{1}{2}\langle\hat x^3\rangle\\
&&\hphantom{{}  +\,}{}\qquad\quad{}\times
\biggl (\biggl\langle\frac{df}{d\hat x}+\hat x\frac{d^2f}{d\hat x^2}\biggr\rangle
+\frac{1}{2}\langle\hat x^3\rangle
\biggl\langle\hat x\frac{df}{d\hat x}+\hat x^2\frac{d^2f}{d\hat
x^2}\biggr\rangle
-\biggl\langle\hat x^2\frac{df}{d\hat x}\biggr\rangle \biggr) \\
&&\qquad\quad{}-\biggl\langle f+\hat x\frac{df}{d\hat x}\biggr\rangle
-\frac{1}{2}\langle\hat x^3\rangle\biggl\langle\hat x f+\hat x^2\frac
{df}{d\hat x}\biggr\rangle
+\langle\hat x^2 f\rangle\\
&&\qquad{}=%1
\biggl\langle\frac{d^2f}{d\hat x^2}\biggr\rangle
+\langle\hat x^3\rangle
\biggl\langle\hat x \frac{d^2f}{d\hat x^2}\biggr\rangle
+\frac{1}{4}\langle\hat x^3\rangle^2
\biggl\langle\hat x^2\frac{d^2f}{d\hat x^2}\biggr\rangle\\
&&\qquad\quad{}+\frac{1}{2}\langle\hat x^3\rangle\biggl\langle\frac{df}{d\hat x}\biggr\rangle
-\frac{1}{2}(1 -2 \langle\hat x^3\rangle^2+\langle\hat x^4\rangle)
\biggl\langle\hat x \frac{df}{d\hat x}\biggr\rangle
-\langle\hat x^3\rangle\biggl\langle\hat x^2\frac{df}{d\hat x}\biggr\rangle\\
&&\qquad\quad{}-\langle f\rangle-\frac{1}{2}\langle\hat x^3\rangle\langle\hat x f\rangle
+\langle\hat x^2 f\rangle.
\end{eqnarray*}

Because of (\ref{19}) we have for any $k \in\N$,
%
%e66 ###
%e64 #&#
\begin{eqnarray}
\label{e_moment_expansion}
\frac{d^k}{d\hat\xi^k}
\biggl(\frac{1}{s}\frac{d}{d\sigma}\biggr)^2\langle \exp(i \hat\xi\hat x)
\rangle
& =&
\biggl(\frac{1}{s}\frac{d}{d\sigma}\biggr)^2\frac{d^k}{d\hat\xi^k}
\langle\exp(i \hat\xi\hat x) \rangle\nonumber
\\[-8pt]
\\[-8pt]
& =&
i^k\biggl(\frac{1}{s}\frac{d}{d\sigma}\biggr)^2
\langle\hat x^k \exp(i \hat\xi\hat x) \rangle.
\nonumber
\end{eqnarray}
This formula and the normalization (\ref{2}) yield that $(\frac
{1}{s}\frac{d}{d\sigma})^2 \langle\exp(i \hat\xi\hat x) \rangle$
vanishes to second order in $\hat\xi$. More precisely, for $k\in
 \{0,1,2  \}$
%
%e67 ###
%e65 #&#
\begin{equation}
\label{20bis}
  \frac{d^k}{d\hat\xi^k} \bigg |_{\hat\xi= 0 }
\biggl(\frac{1}{s}\frac{d}{d\sigma}\biggr)^2 \langle\exp(i \hat\xi\hat x)
\rangle
 =
i^k\biggl(\frac{1}{s}\frac{d}{d\sigma}\biggr)^2
\langle\hat x^k \rangle= 0.
\end{equation}

Therefore, we consider the third derivative w.r.t. $\hat\xi$ given
by \eqref{e_moment_expansion}. For this purpose we apply the formula
for $(\frac{1}{s}\frac{d}{d\sigma})^2\langle f(\hat x)\rangle$ from
above to the function
\[
f= \hat x^3 \exp(i \hat\xi\hat x).
\]
Using the abbreviation $e:=\exp(i\hat\xi\hat x)$, we obtain
%
%d(fg)=df g+ f dg
%d^2(fg)=d^2f g + 2 df dg + f d^2g
%d^3(fg)=d^3f g + 3 d^2f dg + 3 df d^2 g+ f d^3g$
%
\begin{eqnarray*}
\frac{d^3}{d\hat\xi^3} \biggl(\frac{1}{s}\frac{d}{d\sigma}\biggr)^2\langle
e\rangle & = &i^3\biggl(\frac{1}{s}\frac{d}{d\sigma}\biggr)^2
\langle\hat x^3 e \rangle \\
& =& i^3  \biggl( 6 \mv{ \hat x e } + i 6 \hat\xi\mv{\hat x^2 e} -
\hat\xi^2 \mv{\hat x^3 e} \\
&&\hphantom{i^3  \biggl(}{}   + \mv{\hat x^3}  ( 6 \mv{\hat x^2 e} + i 6 \hat
\xi\mv{x^3 e} - \xi^2 \mv{\hat x^4 e} ) \\
&&\hphantom{i^3  \biggl(}{}   + \frac{1}{4} \mv{x^3}^2  ( 6 \mv{\hat x^3 e} +
i 6 \hat\xi\mv{\hat x^4 e} - \hat\xi^2 \mv{\hat x^5 e}  ) \\
&&\hphantom{i^3  \biggl(}{}   + \frac{1}{2} \mv{\hat x^3 }  ( 3 \mv{\hat x^2
e} + i \hat\xi\mv{\hat x^3 e}  ) \\
&&\hphantom{i^3  \biggl(}{}   - \frac{1}{2}  ( 1 - 2 \mv{\hat x^3}^2 + \mv
{\hat x^4}  )  ( 3 \mv{\hat x^3 e} + i \hat\xi\mv{\hat
x^4 e}  ) \\
&&\hphantom{i^3  \biggl(}{}   - \mv{\hat x^3}  ( 3 \mv{\hat x^4 e } + i \hat
\xi\mv{\hat x^5 e}  ) \\
&&\hspace*{55pt}\hphantom{i^3  \biggl(}{}   - \mv{\hat x^3 e} - \frac{1}{2} \mv{\hat x^3} \mv
{\hat x^4 e} + \mv{\hat x^ 5 e}  \biggr).
\end{eqnarray*}
%
% which yields the estimate
% \begin{eqnarray*}
%  | \frac{d^3}{d\hat\xi^3} (\frac{1}{s}\frac{d}{d\sigma})^2\langle
%e \rangle | & \lesssim\mv{| \hat x |} + \mv{| \hat x |^2 }
% ( \mv{ |\hat x|^3 } + | \hat\xi|^2  ) \\
% &  +\mv{ |\hat x|^3 }  ( 1 + \mv{| \hat x|^3}^2 + \mv{ |
% &  + \mv{| \hat x |^4}  ( 1 + \mv{ |\hat x|^3 } + \mv{ |
% &  + \mv{| \hat x |^5}  ( 1 + |\hat\xi| \mv{ |\hat x|^3 }
%+ |\hat\xi|^2 \mv{ |\hat x|^3 }^2  ) .
% \end{eqnarray*}
From this formula and the moment estimates \eqref{3}, we obtain the estimate
\[
\biggl|\frac{d^3}{d\hat\xi^3}
\biggl(\frac{1}{s}\frac{d}{d\sigma}\biggr)^2\langle e\rangle\biggr|
 \lesssim 1+\hat\xi^2.
\]
In combination with (\ref{20bis}), this estimate yields (\ref{20}).
\end{pf*}

%%%%%%%%%%%%%%%%%%%%%%%%%%%%%%%%%%%%%%%%%%%%%%%%%%%%%%%%%%

\begin{pf*}{Proof of Proposition~\ref{P}}
We focus on (\ref{9}) and (\ref{12}). The intermediate estimate (\ref
{11}) can
be established as (\ref{12}).

We start with (\ref{9}). Fix a $\delta>0$ so small such that
the expansion (\ref{C}) of $\langle\exp(i\hat x\hat\xi)\rangle$
holds for $|\hat\xi|\le\delta$.
We split the integral representation
(\ref{13bis}) accordingly:
%
%e68 ###
%e66 #&#
\begin{eqnarray}
\label{P4}
2\pi g_{K, \sigma}(0)&=&
\int_{\{|\fracd{1}{\sqrt{K}}\hat\xi|\le\delta\}}
\biggl\langle\exp\biggl(i\hat x\frac{1}{\sqrt{K}}\hat\xi\biggr)\biggr\rangle^K \,d\hat\xi
\nonumber
\\[-8pt]
\\[-8pt]
&&{}  +
\int_{\{|\fracd{1}{\sqrt{K}}\hat\xi|>\delta\}}
\biggl\langle\exp\biggl(i\hat x\frac{1}{\sqrt{K}}\hat\xi\biggr)\biggr\rangle^K \,d\hat\xi
=: I + \mathit{II}.
\nonumber
\end{eqnarray}

We consider the first term $I$ on the r.h.s. of (\ref{P4}),
which will turn out to be of leading order. Since
$\delta$ is so small that (\ref{C}) holds, we may rewrite it as
%
%e69 ###
%e67 #&#
\begin{eqnarray}\label{P6}
I&:=&\int_{\{|\fracd{1}{\sqrt{K}}\hat\xi|\le\delta\}}
\biggl\langle\exp\biggl(i\hat x\frac{1}{\sqrt{K}}\hat\xi\biggr)\biggr\rangle^K
d\hat\xi\nonumber
\\[-8pt]
\\[-8pt]
&\hspace*{2.8pt}=&
\int_{\{|\fracd{1}{\sqrt{K}}\hat\xi|\le\delta\}}
\exp\biggl(-Kh\biggl(\frac{1}{\sqrt{K}}\hat\xi\biggr)\biggr)\,d\hat\xi.
\nonumber
\end{eqnarray}
We note that for $|\frac{1}{\sqrt{K}}\hat\xi|\le\delta$ we have
by (\ref{C}),
%
%e70 ###
%e68 #&#
\begin{equation}\label{P5}
\biggl|Kh\biggl(\frac{1}{\sqrt{K}}\hat\xi\biggr)-\frac{1}{2}\hat\xi^2\biggr|
 \lesssim \frac{1}{\sqrt{K}}|\hat\xi|^3,
\end{equation}
in particular for $\delta$ small enough,
%
%e71 ###
%e69 #&#
\begin{equation}\label{P2}
\Real\biggl ( Kh\biggl(\frac{1}{\sqrt{K}}\hat\xi\biggr)  \biggr)  \ge \frac
{1}{4}\hat\xi^2,
\end{equation}
so that (\ref{P5}) implies by the Lipschitz continuity of
$\mathbb{C} \ni y\mapsto\exp(y) \in\mathbb{C}$ on $ \Real y\le
-\frac{1}{4}\hat\xi^2$ with constant
$\exp(-\frac{1}{4}\hat\xi^2)$,
\[
\biggl|\exp\biggl(-Kh\biggl(\frac{1}{\sqrt{K}}\hat\xi\biggr)\biggr)-\exp\biggl(-\frac{1}{2}\hat\xi^2\biggr)\biggr|
 \lesssim \frac{1}{\sqrt{K}}|\hat\xi|^3\exp\biggl(-\frac{1}{4}\hat
\xi^2\biggr).
\]
Inserting this estimate into (\ref{P6}) we obtain
\begin{eqnarray*}
\biggl|I-\int_{\{|\fracd{1}{\sqrt{K}}\hat\xi|\le\delta\}}
\exp\biggl(-\frac{1}{2}\hat\xi^2\biggr)\,d\hat\xi\biggr|
  & \lesssim&
\frac{1}{\sqrt{K}}\int_{\{|\fracd{1}{\sqrt{K}}\hat\xi|\le\delta
\}}
|\hat\xi|^3\exp\biggl(-\frac{1}{4}\hat\xi^2\biggr)\,d\hat\xi\\
& \lesssim&  \frac{1}{\sqrt{K}}\int
|\hat\xi|^3\exp\biggl(-\frac{1}{4}\hat\xi^2\biggr)\,d\hat\xi\\
& \lesssim&  \frac{1}{\sqrt K}.
\end{eqnarray*}
The latter turns, as desired, into
\begin{eqnarray}
\bigl|I-\sqrt{2\pi}\bigr|
 &=&
\biggl|I-\int\exp\biggl(-\frac{1}{2}\hat\xi^2\biggr)\,d\hat\xi\biggr|\nonumber\\
&\lesssim&
\frac{1}{\sqrt{K}}
+\int_{\{|\fracd{1}{\sqrt{K}}\hat\xi|>\delta\}}
\exp\biggl(-\frac{1}{2}\hat\xi^2\biggr)\,d\hat\xi\nonumber\\
&\lesssim&\frac{1}{\sqrt{K}}, \nonumber
\end{eqnarray}
since $\int_{\{|\fracd{1}{\sqrt{K}}\hat\xi|>\delta\}}
\exp(-\frac{1}{2}\hat\xi^2)\,d\hat\xi$ is exponentially small in $K$.

We now address the second term $\mathit{II}$ on the r.h.s. of (\ref{P4}); on
the integrand we use Lemma~\ref{soft} (on $K-2$ of the $K$ factors)
and (\ref{L2}) (on the remaining $2$ factors).
\begin{eqnarray}
\biggl|\biggl\langle\exp\biggl(i\hat x\frac{1}{\sqrt{K}}\hat\xi\biggr)\biggr\rangle\biggr|^{K}
&\lesssim&
\lambda^{K-2}   \biggl(\frac{1}{1+ \fracd{1}{\sqrt{K}} |\xi|}
 \biggr)^{2}
\nonumber\\
& \lesssim&
K  \lambda^{K-2}   \frac{1}{K+\hat\xi^{2}}   \lesssim
K  \lambda^{K-2}  \frac{1}{1+\hat\xi^2}.\nonumber
\end{eqnarray}
It follows that the second term $\mathit{II}$ on the r.h.s. of (\ref{P4}) is
exponentially small and thus of higher order:
\begin{eqnarray}
\biggl| \int_{\{|\fracd{1}{\sqrt{K}}\hat\xi|>\delta\}}
\biggl\langle\exp\biggl(i\hat x\frac{1}{\sqrt{K}}\hat\xi\biggr)\biggr\rangle^K d\hat\xi
\biggr|&\lesssim& K \lambda^{K-2}  \int\frac{1}{1+\hat\xi^2}\,d\hat\xi\nonumber\\
&\lesssim& K \lambda^{K-2}  \stackrel{\lambda<1}{\ll} \frac{1}{\sqrt{K}}.\nonumber
\end{eqnarray}

We now turn to (\ref{12}). We take the second
$\sigma$-derivative of the integral representation
(\ref{13bis}),
\begin{eqnarray}
&&2\pi\biggl(\frac{1}{s}\frac{d}{d\sigma}\biggr)^2g_{K, \sigma
}(0)\nonumber\\
&& \qquad =
\int\biggl (K(K-1)\biggl\langle\exp\biggl(i\hat x\frac{1}{\sqrt{K}}\hat\xi
\biggr)\biggr\rangle^{K-2}
\biggl(\frac{1}{s}\frac{d}{d\sigma}\biggl\langle
\exp\biggl(i\hat x\frac{1}{\sqrt{K}}\hat\xi\biggr)\biggr\rangle\biggr)^2\nonumber\\
&& \qquad  \quad {}+K\biggl\langle\exp\biggl(i\hat x\frac{1}{\sqrt{K}}\hat\xi\biggr)\biggr\rangle^{K-1}
\biggl(\frac{1}{s}\frac{d}{d\sigma}\biggr)^2\biggl\langle\exp\biggl(i\hat x\frac{1}{\sqrt
{K}}\hat\xi
\biggr)\biggr\rangle \biggr)\,d\hat\xi\nonumber
\end{eqnarray}
and use Lemma~\ref{L3},
%
%e72 ###
%e70 #&#
\begin{eqnarray}\label{P1}
   \biggl|\biggl(\frac{1}{s}\frac{d}{d\sigma}\biggr)^2g_{K, \sigma}(0) \biggr|   &\lesssim&
\int \biggl(K^2\biggl|\biggl\langle\exp\biggl(i\hat x\frac{1}{\sqrt{K}}\hat\xi
\biggr)\biggr\rangle\biggr|^{K-2}
\biggl(1+\biggl|\frac{1}{\sqrt{K}}\hat\xi\biggr|^2\biggr)\biggl|\frac{1}{\sqrt{K}}\hat\xi\biggr|^6
\nonumber\\
&&      {}     +K\biggl|\biggl\langle\exp\biggl(i\hat x\frac{1}{\sqrt{K}}\hat\xi
\biggr)\biggr\rangle\biggr|^{K-1}
\biggl(1+\biggl|\frac{1}{\sqrt{K}}\hat\xi\biggr|^2\biggr)\biggl|\frac{1}{\sqrt{K}}\hat\xi
\biggr|^3 \biggr)\,d\hat\xi
 \\
    &\lesssim&    \frac{1}{\sqrt{K}}
\int\biggl|\biggl\langle\exp\biggl(i\hat x\frac{1}{\sqrt{K}}\hat\xi\biggr)\biggr\rangle\biggr|^{K-2}
\biggl(1+\biggl|\frac{1}{\sqrt{K}}\hat\xi\biggr|^2\biggr)(|\hat\xi|^6+1)\,d\hat\xi.
\nonumber
\end{eqnarray}
As for (\ref{9}), we split the integral representation (\ref{P1})
according to $\delta$:
%
%e73 ###
%e71 #&#
\begin{eqnarray}\label{P3}
&&  \biggl|\biggl(\frac{1}{s}\frac{d}{d\sigma}\biggr)^2g_{K, \sigma}(0) \biggr| \nonumber\hspace*{-30pt}\\
&& \quad   \lesssim  \frac{1}{\sqrt{K}}
\int_{\{\fracd{1}{\sqrt{K}}|\hat\xi|\le\delta\}}
\biggl|\biggl\langle\exp\biggl(i\hat x\frac{1}{\sqrt{K}}\hat\xi\biggr)\biggr\rangle\biggr|^{K-2}
\biggl(1+\biggl|\frac{1}{\sqrt{K}}\hat\xi\biggr|^2\biggr)(\hat\xi^6+1)\,d\hat\xi\nonumber\hspace*{-30pt}\\
&& \qquad{}      + \frac{1}{\sqrt{K}}
\int_{\{\fracd{1}{\sqrt{K}}|\hat\xi|>\delta\}}
\biggl|\biggl\langle\exp\biggl(i\hat x\frac{1}{\sqrt{K}}\hat\xi\biggr)\biggr\rangle\biggr|^{K-2}
\biggl(1+\biggl|\frac{1}{\sqrt{K}}\hat\xi\biggr|^2\biggr)(\hat\xi^6+1)\,d\hat\xi\hspace*{-30pt}\\
&& \quad \lesssim  \frac{1}{\sqrt{K}}
\int_{\{\fracd{1}{\sqrt{K}}|\hat\xi|\le\delta\}}
\biggl|\biggl\langle\exp\biggl(i\hat x\frac{1}{\sqrt{K}}\hat\xi\biggr)\biggr\rangle\biggr|^{K-2}
(\hat\xi^6+1)\,d\hat\xi\nonumber\hspace*{-30pt}\\
&& \qquad{}     +\frac{1}{\sqrt{K}}
\int_{\{\fracd{1}{\sqrt{K}}|\hat\xi|>\delta\}}
\biggl|\biggl\langle\exp\biggl(i\hat x\frac{1}{\sqrt{K}}\hat\xi\biggr)\biggr\rangle\biggr|^{K-2}
(\hat\xi^8+1)\,d\hat\xi.\nonumber\hspace*{-30pt}
\end{eqnarray}
On the first r.h.s. term we use (\ref{P2}):
\begin{eqnarray}
&&\frac{1}{\sqrt{K}}\int_{\{\fracd{1}{\sqrt{K}}|\hat\xi
|\le\delta\}}
\biggl|\biggl\langle\exp\biggl(i\hat x\frac{1}{\sqrt{K}}\hat\xi\biggr)\biggr\rangle\biggr|^{K-2}
(\hat\xi^6+1)\,d\hat\xi\nonumber\\
&& \qquad \lesssim
\frac{1}{\sqrt{K}}\int_{\{\fracd{1}{\sqrt{K}}|\hat\xi|\le\delta
\}}
\exp\biggl(-(K-2)\frac{1}{4}\biggl(\frac{1}{\sqrt{K}}\hat\xi\biggr)^2\biggr)
(\hat\xi^6+1)\,d\hat\xi\nonumber\\
&& \qquad \stackrel{K\gg1}{\lesssim}
\frac{1}{\sqrt{K}}\int\exp\biggl(-\frac{1}{8}\hat\xi^2\biggr)(\hat\xi
^6+1)\,d\hat\xi
\nonumber\\
&& \qquad \lesssim\frac{1}{\sqrt{K}}. \nonumber
\end{eqnarray}
On the integrand of the second r.h.s. term in (\ref{P3})
we use Lemma~\ref{soft} (on $K-12$ of the $K-2$ factors) and (\ref{L2})
(on the remaining $10$ factors):
\begin{eqnarray}
\biggl|\biggl\langle\exp\biggl(i\hat x\frac{1}{\sqrt{K}}\hat\xi\biggr)\biggr\rangle\biggr|^{K-2}
(\hat\xi^8+1)
&\lesssim&
\lambda^{K-12}   \biggl(\frac{1}{1+ \fracd{1}{\sqrt{K}} |\xi|}
 \biggr)^{10}(\hat\xi^8+1)
\nonumber\\
& \lesssim&
K^{5}\lambda^{K-12}   \frac{1}{K^{5}+\hat\xi^{10}}  (\hat\xi
^8+1)\nonumber\\
&\lesssim&
K^{5}\lambda^{K-12}\frac{1}{1+\hat\xi^2}.\nonumber
\end{eqnarray}
Hence, we see that this second term in (\ref{P3}) is exponentially small
and thus of higher order:
\begin{eqnarray}
&&\frac{1}{\sqrt{K}}
\int_{\{\fracd{1}{\sqrt{K}}|\hat\xi|>\delta\}}
\biggl|\biggl\langle\exp\biggl(i\hat x\frac{1}{\sqrt{K}}\hat\xi\biggr)\biggr\rangle\biggr|^{K-2}
(|\hat\xi|^8+1)\,d\hat\xi\nonumber\\
&& \qquad \lesssim
K^{9/2}\lambda^{K-12}\int\frac{1}{1+\hat\xi^2}\,d\hat\xi\nonumber
\\
&& \qquad  \lesssim K^{9/2}\lambda^{K-12} \stackrel{\lambda<1}{\ll} \frac
{1}{\sqrt{K}}. \nonumber
\end{eqnarray}
\upqed
\end{pf*}

%%%%%%%%%%%%%%%%%%%%%%%%%%%%%%%%%%%%%%%%%%%%%%%%%%%%%%%%%%%%%%%%%%%%%%

For the proof of Lemma~\ref{p_assumptions_lct} we need the following
auxiliary statement, based on elementary calculus.
%
%le3.6 #&#
\begin{lemma}\label{p_aux_char}
Assume that the single-site potential $\psi\dvtx \R\to\R$ is convex. We
consider the corresponding Gibbs measure,
\[
\nu(dx) = \frac{1}{Z}  \exp(- \psi(x))\,dx.
\]
Let $M$ denote the maximum of the density of $\nu$, that is,
\[
M := \max_{x} \frac{1}{Z} \exp(- \psi(x)).
\]
Then we have for all $k \in\N$,
\begin{eqnarray*}
\int|x|^k  \nu(dx)  \lesssim \frac{1}{M^k}
\end{eqnarray*}
for some constant only depending on $k$.
\end{lemma}

\begin{pf*}{Proof of Lemma~\ref{p_aux_char}}
We may assume w.l.o.g. that
%
%e74 ###
%e72 #&#
\begin{equation}\label{e_mass_finite}
Z= \int\exp(-\psi(x))\,dx=1,
\end{equation}
and $M:= \sup_x \exp(- \psi(x))$ is attained at $x= 0$, which means
%
%e75 ###
%e73 #&#
\begin{equation} \label{e_d_M}
M= \exp(- \psi(0)).
\end{equation}
It follows from convexity of $\psi$ that
%
%e76 ###
%e74 #&#
\begin{eqnarray}\label{e_cos_conv}
\psi'(x) \leq0   \qquad \mbox{for }  x \leq0 \quad\mbox{and} \quad
\psi'(x) \geq0    \qquad \mbox{for }   x \geq0.
\end{eqnarray}
We start with an analysis of the convex single-site potential $\psi$.
We first argue that
%
%e77 ###
%e75 #&#
\begin{equation} \label{e_convex_inner_bound}
\psi \biggl( \pm\frac{e}{M}  \biggr) \geq- \log M + \log e = - \log
M + 1.
\end{equation}
Indeed in view of the monotonicity \eqref{e_cos_conv}, we have
\begin{eqnarray*}
1 \overset{\mathclap{\eqref{e_mass_finite}}}{\geq} \int_0^{\fraca{e}{M}} \exp
(- \psi(y))\,dy \overset{\mathclap{\eqref{e_cos_conv}}}{\geq} \frac{e}{M}
\exp\biggl ( - \psi \biggl( \frac{e}{M}  \biggr) \biggr)
\end{eqnarray*}
and
\begin{eqnarray*}
1 \overset{\mathclap{\eqref{e_mass_finite}}}{\geq} \int_{-\fraca{e}{M}}^0 \exp
(- \psi(y))\,dy \overset{\mathclap{\eqref{e_cos_conv}}}{\geq} \frac{e}{M}
\exp\biggl ( - \psi \biggl(- \frac{e}{M}  \biggr) \biggr).
\end{eqnarray*}
We now argue that for $|x| \geq\frac{e}{M}$,
%
%e78 ###
%e76 #&#
\begin{equation} \label{e_lower_bound_psi}
\psi(x)\geq\frac{M}{e}  \biggl(|x| - \frac{e}{M}  \biggr) - \log M.
\end{equation}
W.l.o.g. we may restrict ourselves to $x \geq\frac{e}{M}$. By
convexity of $\psi$, we have
\[
\psi'  \biggl( \frac{e}{M}  \biggr)  \frac{e}{M} \geq \psi \biggl(
\frac{e}{M}  \biggr) - \psi(0) \overset{\mathclap{\eqref
{e_d_M}}}{=}  \psi\biggl ( \frac{e}{M}  \biggr) + \log M \overset
{\mathclap{\eqref{e_convex_inner_bound}}}{\geq}  1.
\]
The convexity of $\psi$, the last estimate and \eqref
{e_convex_inner_bound} yield for $ x \geq\frac{e}{M}$, as desired,
\begin{eqnarray*}
\psi(x) & \geq&\psi' \biggl ( \frac{e}{M}  \biggr)  \biggl(x - \frac
{e}{M}  \biggr) + \psi\biggl ( \frac{e}{M}  \biggr) \\
& \geq&\frac{M}{e}  \biggl(x - \frac{e}{M}  \biggr) - \log M.
\end{eqnarray*}

We finished the analysis on $\psi$ and turn to the verification of the
estimate of Lemma~\ref{p_aux_char}. We split the integral according to
\[
\int|x|^k \exp(- \psi(x))\,dx = \int_{-\infty}^0 |x|^k \exp(- \psi
(x))\,dx + \int_0^\infty|x|^k \exp(- \psi(x))\,dx.
\]
We will now deduce the estimate
\[
\int_0^\infty|x|^k \exp(- \psi(x))\,dx \lesssim\frac{1}{M^k}.
\]
A similar estimate for the integral $\int_{-\infty}^0 |x|^k \exp(-
\psi(x))\,dx$ follows from the same argument by symmetry. We split the integral
\begin{eqnarray*}
&&\int_0^\infty |x|^k \exp(- \psi(x))\,dx\\
&& \qquad  = \int_0^{\fraca{e}{M}}
|x|^k \exp(- \psi(x))\,dx + \int_{\fraca{e}{M}}^\infty |x|^k \exp(-
\psi(x))\,dx.
\end{eqnarray*}
The first integral on the r.h.s. can be estimated as
\[
\int_0^{\fraca{e}{M}} |x|^k \exp(- \psi(x))\,dx \leq \frac{e^k}{M^k}
\int\exp(- \psi(x))\,dx \overset{\mathclap{\eqref{e_mass_finite}}}{=} \frac{e^k}{M^k}.
\]
For the estimation of the second integral, we apply \eqref
{e_lower_bound_psi}, which yields, by the change of variables $\frac
{M}{e}  ( x - \frac{e}{M} ) = \hat x$,
\begin{eqnarray*}
\int_{\fraca{e}{M}}^\infty |x|^k \exp(- \psi(x))\,dx &\leq& \int
_{\fraca{e}{M}}^\infty |x|^k \exp \biggl(- \frac{M}{e}  \biggl(x -
\frac{e}{M}  \biggr) + \log M  \biggr)\,dx \\
& =& M  \frac{e}{M} \int_{0}^\infty  \biggl| \frac{e}{M} \hat x
+\frac{e}{M}  \biggr|^k \exp (- \hat x  )\,d\hat x \\
& =& e  \biggl ( \frac{e}{M}  \biggr)^k \int_{0}^\infty | \hat x + 1
|^k \exp (- \hat x  )\,d\hat x \\
& \lesssim&\frac{1}{M^k}.
\end{eqnarray*}
\upqed
\end{pf*}

Equipped with Lemma~\ref{p_aux_char}, we are able to give an
elementary proof of Lemma~\ref{p_assumptions_lct}:

\begin{pf*}{Proof of Lemma~\ref{p_assumptions_lct}}
We argue that $s \lesssim1$. Because $\psi$ is a bounded perturbation
of a uniformly strictly convex function, the measure $\mu^{\sigma}$
given by \eqref{e_definition_mu_sigma} satisfies the SG uniformly in
$\sigma$. This implies, in particular,
%
%e79 ###
%e77 #&#
\begin{equation}\label{e_est_s}
s^2 = \var_{\mu^{\sigma}} (x) \lesssim\int \biggl( \frac{d}{dx}
x  \biggr)^2\,d\mu^{\sigma} = 1
\end{equation}
uniformly in $\sigma$ and thus in $m$.

Now, we verify \eqref{4}. Using $|\delta\psi|\lesssim1$ to pass
from $\psi$ to $\psi_c$,
we may assume that $\psi$ is strictly convex. In fact, we can give up
\emph{strict} convexity of $\psi$ and may only assume that $\psi$ is
convex. % We assume w.l.o.g. that
% \[
% \int exp(- \psi(x)) = 1.
% \]
By the change of variables $\hat x = \frac{x-m}{s}$, we have for any $
k \in\N$,
\[
\frac{ \int|x-m|^k\,d\mu}{s^k} = \int|\hat x |^k \exp(- \hat\psi
(\hat x))\,d\hat x
\]
for some convex function $\hat\psi$, which is normalized in the sense that
%
%e80 ###
%e78 #&#
\begin{equation}
\int\exp(- \hat\psi(\hat x))\,d\hat x = 1  \quad\mbox{and}
\quad  \int\hat{x}^2 \exp(- \hat\psi(\hat x))\,d\hat x = 1.
\label{e_norm_psi}
\end{equation}
An application of Lemma~\ref{p_aux_char} yields the estimate
\[
\frac{ \int|x-m|^k\,d\mu}{s^k} \leq\int|\hat x |^k \exp(- \hat\psi
(\hat x))\,d\hat x \lesssim\frac{1}{M^k},
\]
where $M$ is given by $M:= \max_{\hat x} \exp(- \hat\psi(\hat x))$.
Now, we argue that due to the normalization of $\hat\psi$, we have
\[
M \geq C
\]
for some universal constant $C>0$. The latter verifies the desired
estimate~\eqref{4}. Indeed normalization \eqref{e_norm_psi} implies
\begin{eqnarray}
\int_{(-2,2)}\exp(-\psi(\hat x))\,d\hat x&\stackrel{\mathclap{\eqref{e_norm_psi}}}{=}&
1-\int_{\mathbb{R}-(-2,2)}\exp(-\psi(\hat x))\,d\hat x\nonumber\\
&\ge&
1-\frac{1}{4}\int\hat x^2\exp(-\psi(\hat x))\,d\hat x
 \stackrel{\mathclap{\eqref{e_norm_psi}}}{\ge} \frac{3}{4}.\nonumber
\end{eqnarray}
Hence, there exists an $\hat x_0\in(-2,2)$ such that
$\exp(-\psi(\hat x_0))\ge\frac{3}{8}$, which yields
\[
M= \max_{\hat x} \exp(- \hat\psi(\hat x)) \geq\exp(-\psi(\hat
x_0))\ge\frac{3}{8}.
\]
Let us turn to the statement \eqref{7} of Proposition~\ref{P}. Writing
\[
\exp ( i x \xi ) = \frac{d}{dx} \biggl ( -i  \frac
{1}{\xi}  \exp(ix \xi) \biggr),
\]
we obtain by integration by parts that
\begin{eqnarray*}
\mv{\exp ( i x \xi ) } & =& i  \frac{1}{\xi} \int\exp
 ( i x \xi ) \, \frac{d}{dx}  \bigl( \exp\bigl ( - \varphi
^*(\sigma) + \sigma x - \psi(x)  \bigr)  \bigr)\,dx \\
& =& i  \frac{1}{\xi} \int\exp ( i x \xi )    \bigl(
\sigma - \psi' (x)  \bigr) \exp\bigl ( - \varphi^*(\sigma) +
\sigma x - \psi(x)  \bigr)\,dx .
\end{eqnarray*}
For convenience, we introduce the Hamiltonian $\hat\psi(x) =-\sigma x
+ \psi_c (x)$ and assume w.l.o.g. that $\int\exp(- \hat\psi(x))\,dx
= 1$.
The splitting $\psi= \psi_c + \delta\psi$ with $| \delta\psi |$,
$|\delta\psi'| \lesssim1 $ and definition \eqref
{e_definition_varphi*} of $\varphi^*$ yield the estimate
\begin{eqnarray*}
 | \mv{\exp ( i x \xi ) }  | & \lesssim &
\frac{1}{|\xi|}   \frac{  \int  | \sigma - \psi_c' (x)-
\delta\psi_c' (x)  | \exp ( \sigma x - \psi_c (x)-
\delta\psi_c (x)  )\,dx} { \int \exp ( \sigma x - \psi_c
(x) - \delta\psi_c (x)  )\,dx } \\
& \lesssim&  \frac{1}{s|\xi|}   s  \int | \hat\psi' (x) | \exp
( - \hat\psi(x) )\,dx + \frac{1}{s | \xi|}  s ,
\end{eqnarray*}
where $s$ is defined as in Proposition~\ref{P}. Because $s \lesssim1$
by \eqref{e_est_s}, we only have to consider the first term of the
r.h.s. of the last inequality. We argue that for
\[
M := \max_x \exp( - \hat\psi(x) ),
\]
it holds
%
%e81 ###
%e79 #&#
\begin{equation}
2 M = \int | \hat\psi' (x) | \exp( - \hat\psi(x) )\,dx. \label{e_f_M}
\end{equation}
For the proof of the last statement, we only need the fact that $\hat
\psi(x) =-\sigma x + \psi_c (x)$ is convex. W.l.o.g. we may assume
that $M$ is attained at $x= 0$, which means $M= \exp(- \hat\psi
(0))$. It follows from convexity of $\hat\psi$ that
\begin{eqnarray*}
\hat\psi'(x) \leq0   \qquad \mbox{for }  x \leq0 \quad\mbox{and}
\quad
\hat\psi'(x) \geq0   \qquad  \mbox{for }   x \geq0. % \label{e_cos_conv}
\end{eqnarray*}
Indeed, we get
\begin{eqnarray*}
&&\int|\hat\psi'(x)|  \exp(- \hat\psi(x))\,dx \\
&& \qquad  = - \int_{- \infty}^0 \hat\psi'(x) \exp(- \hat\psi(x))\,dx +
\int_0^\infty\hat\psi'(x) \exp( - \hat\psi(x))\,dx \\
&& \qquad  = 2 \exp(- \hat\psi(0)) = 2 M.
\end{eqnarray*}
Because the mean of a measure $\mu$ is optimal in the sense that for
all $c \in\R$,
\[
% \label{e_mean_optimal}
%
\int (x-c  )^2 \mu(dx) \geq\int \biggl( x - \int x \mu
(dx)  \biggr)^2 \mu(dx),
\]
we can estimate
%
%e82 ###
%e80 #&#
\begin{eqnarray}\label{e_var_2nd_mom}
s^2 \leq\frac{\int x^2 \exp ( \sigma x - \psi(x)  )\,dx}{\int \exp ( \sigma x - \psi(x)  )\,dx}  \overset
{|\delta\psi| \lesssim1 }{\lesssim}  \int x^2 \exp( - \hat\psi
(x) )\,dx .
\end{eqnarray}
Therefore, Lemma~\ref{p_aux_char} applied to $k=2$ and $\psi$
replaced by $\hat\psi$ yields
\[
s  \int  | \hat\psi' (x)  | \exp ( - \hat\psi(x)
 )\,dx \overset{\mathclap{\eqref{e_f_M} ,\eqref{e_var_2nd_mom}}}{\lesssim
}   \biggl ( \int x^2 \exp( - \hat\psi(x) )\,dx  \biggr)^\h M
\lesssim  1,
\]
which verifies \eqref{7} of Proposition~\ref{P}.
\end{pf*}

Before we turn to the proof of Lemma~\ref{p_crucial_convexification},
we will deduce the following auxiliary result.
%
%le3.7 #&#
\begin{lemma}\label{e_aux_est_d_dm_s}
Assume that \eqref{4} of Proposition~\ref{P} is satisfied. Then,
using the notation of Proposition~\ref{P}, it holds that
\[
\textup{(i)} \quad    \biggl| \frac{d}{dm} s  \biggr| \lesssim1  \quad
\mbox{and}   \quad\textup{(ii)} \quad    \biggl| \frac{d^2}{dm^2} s
 \biggr| \lesssim \frac{1}{s}.
\]
\end{lemma}

\begin{pf*}{Proof of Lemma~\ref{e_aux_est_d_dm_s}} We start with
restating some basic identities [cf. \eqref{e_d_dsigma_m} and \eqref
{e_d^2_d^2sigma_m}]: It holds that
%
%e85 ###
%e84 ###
%e83 ###
%e81 #&#
%e82 #&#
%e83 #&#
\begin{eqnarray}\label{e_d_dsigma_m_2}
\frac{d}{d \sigma} m &=& s^2,  \\\label{e_d^2_dsigma^2_m}
\frac{d^2}{d \sigma^2 } m & =& \frac{d}{d \sigma} s^2 = \int (
x- m  )^3 \mu^\sigma(dx),  \\\label{e_d^3_dsigma^3_m}
\frac{d^3}{d \sigma^3 } m &=& \int ( x- m  )^4 \mu^\sigma
(dx).
\end{eqnarray}

Let us consider (i): It follows from \eqref{e_d_dsigma_m_2}
and \eqref{e_d^2_dsigma^2_m} that
\begin{eqnarray*}
\frac{d}{dm} s^2 & = &\frac{d}{d \sigma} s^2 \,\frac{d}{dm} \sigma\\
& = &\int ( x- m  )^3 \mu^\sigma(dx)  \biggl(\frac{d}{d
\sigma} m  \biggr)^{-1} \\
& =& \frac{\int ( x- m  )^3 \mu^\sigma(dx)}{s^3 }  s,
\end{eqnarray*}
which yields by assumption \eqref{4} of Proposition~\ref{P} the estimate
\[
 \biggl| \frac{d}{dm} s^2  \biggr| \lesssim s.
\]
The statement of (i) is a direct consequence of the last estimate and
the identity
\[
\frac{d}{dm} s = \frac{1}{2s} \frac{d}{dm} s^2.
\]

We turn to statement (ii): Differentiating the last identity yields
\[
\frac{d^2}{dm^2} s = - \frac{1}{2} \frac{1}{s^2} \frac{d}{dm} s
\frac{d}{dm} s^2 + \frac{1}{2s} \frac{d^2}{dm^2} s^2.
\]
The estimation of the first term on the r.h.s. follows from the estimates
\[
 \biggl| \frac{d}{dm} s^2  \biggr| \lesssim s \quad\mbox{and}
\quad  \biggl| \frac{d}{dm} s  \biggr| \lesssim1,
\]
which we have deduced in the first step of the proof. We turn to the
estimation of the second term. A direct calculation using \eqref
{e_d_dsigma_m_2} yields the identity
%
%e86 ###
%e84 #&#
\begin{eqnarray}
\label{e_i_d^2_dm^2_s^2}
\frac{d^2}{dm^2} s^2 &=& \frac{d^2}{dm^2} \frac{d}{d\sigma} m  =
\frac{d}{dm} \biggl ( \frac{d^2}{d \sigma^2} m \frac{d}{dm} \sigma
 \biggr)\nonumber
 \\[-8pt]
 \\[-8pt]
& =& \frac{d^3}{d \sigma^3} m  \biggl( \frac{d}{dm} \sigma \biggr)^2
+ \frac{d^2}{d \sigma^2} m \frac{d^2}{dm^2} \sigma.
\nonumber
\end{eqnarray}
Considering the first term on the r.h.s., we get from the
identities \eqref{e_d_dsigma_m_2} and~\eqref{e_d^3_dsigma^3_m}, and
the assumption \eqref{4} of Proposition~\ref{P} that
\[
 \biggl| \frac{d^3}{d \sigma^3} m \biggl ( \frac{d}{dm} \sigma
\biggr)^2  \biggr| = \frac{\int ( x- m  )^4 \mu^\sigma
(dx)}{s^4} \lesssim1.
\]
Before we consider the second term of the r.h.s. of \eqref
{e_i_d^2_dm^2_s^2}, we establish the following estimate:
%
%e87 ###
%e85 #&#
\begin{equation}
\label{e_d^2_dm^2_sigma}
 \biggl| \frac{d^2}{dm^2} \sigma  \biggr| \lesssim\frac{1}{s^3}.
\end{equation}
Indeed, direct calculation using \eqref{e_d_dsigma_m_2} and \eqref
{e_d^2_dsigma^2_m} yields
\begin{eqnarray*}
\frac{d^2}{dm^2} \sigma& = & \biggl( \frac{d}{d\sigma} \frac{d}{dm}
\sigma \biggr) \,\frac{d}{dm} \sigma\\
& =&  \biggl( \frac{d}{d\sigma} \biggl ( \frac{d}{d\sigma} m
\biggr)^{-1} \biggr)  \biggl( \frac{d}{d \sigma} m  \biggr)^{-1} \\
& =& - \biggl ( \frac{d}{d \sigma} m  \biggr)^{-3} \frac{d^2}{d \sigma
^2} m \\
& =& - \frac{1}{s^3} \frac{\int ( x- m  )^3 \mu^\sigma(dx)}{s^3}.
\end{eqnarray*}
The last identity yields \eqref{e_d^2_dm^2_sigma} using the
assumption \eqref{4} of Proposition~\ref{P}. Using \eqref
{e_d^2_dm^2_sigma} and \eqref{e_d^2_dsigma^2_m}, we can estimate the
second term of the r.h.s. of \eqref{e_i_d^2_dm^2_s^2} as
\begin{eqnarray*}
 \biggl| \frac{d^2}{d \sigma^2} m \,\frac{d^2}{dm^2} \sigma \biggr|
\lesssim\frac{1}{s^3}   \biggl| \int ( x- m  )^3 \mu
^\sigma(dx)  \biggr|.
\end{eqnarray*}
By applying assumption \eqref{4} of Proposition~\ref{P} this yields
\[
 \biggl| \frac{d^2}{d \sigma^2} m \,\frac{d^2}{dm^2} \sigma \biggr|
\lesssim 1,
\]
which concludes the argument for (ii).
\end{pf*}

\begin{pf*}{Proof of Lemma~\ref{p_crucial_convexification}}
Recall the representation \eqref{e_cramer_idea}, that is,
\[
\tilde g_{K,m}(0)= \exp \bigl(K \varphi(m) - K \bar H_K (m)  \bigr) .
\]
Here, $\tilde g_{K,m}(\xi)$ denotes the Lebesgue density of the random variable
\[
\frac{1}{\sqrt{K}} \sum_{i=1}^K  ( X_i - m  ),
\]
where $X_i$ are real-valued independent random variables identically
distributed according to $\mu^{\sigma}$; cf. \eqref
{e_definition_mu_sigma}. Let $g_{K, \sigma}$ denote the density of the
normalized random variable
\[
\frac{1}{\sqrt{K}} \sum_{i=1}^K \frac{X_i - m}{s} ,
\]
where $s$ is given by \eqref{1}. Then the densities are related by
\[
\frac{1}{s}g_{K, \sigma}  \biggl( \frac{x}{s}  \biggr) = \tilde
g_{K,m} (x).
\]
It follows from \eqref{e_cramer_idea} that
\[
K \varphi(m) - K \bar H_K (m) = \log g_{K, \sigma}(0) - \log s.
\]
In order to deduce the desired estimate, it thus suffices to show
%
%e88 ###
%e86 #&#
\begin{equation}
\label{e_d^2_dm_log_s}
 \biggl| \frac{d^2}{dm^2} \log s  \biggr| \lesssim\frac{1}{s^2}
\end{equation}
and
%
%e89 ###
%e87 #&#
\begin{equation}
\label{e_d^2_dm_log_g}
 \biggl| \frac{d^2}{dm^2} \log g_{K, \sigma} (0)  \biggr| \lesssim
\frac{1}{s^2}.
\end{equation}

The first estimate follows directly from the identity
\begin{eqnarray*}
\frac{d^2}{dm^2} \log s = \frac{d}{dm} \biggl ( \frac{1}{s} \frac
{d}{dm} s  \biggr) = - \frac{1}{s^2} \biggl ( \frac{d}{dm} s
\biggr)^2 + \frac{1}{s} \frac{d^2}{dm^2} s
\end{eqnarray*}
and the estimates provided by Lemma~\ref{e_aux_est_d_dm_s}.

We turn to the second estimate. The identity
\[
\frac{d^2}{dm^2} \log g_{K, \sigma} = - \frac{1}{g_{K, \sigma}^2}
 \biggl( \frac{d}{dm} g_{K, \sigma}  \biggr)^2 + \frac{1}{g_{K,
\sigma}} \frac{d^2}{dm^2} g_{K, \sigma}
\]
and \eqref{9} yield for large $K$ the estimate
\[
 \biggl| \frac{d^2}{dm^2} \log g_{K, \sigma} (0)  \biggr| \lesssim
\biggl ( \frac{d}{dm} g_{K, \sigma} (0)  \biggr)^2 +  \biggl| \frac
{d^2}{dm^2} g_{K, \sigma} (0)  \biggr|.
\]
The estimation of the first term on the r.h.s. follows from
estimate \eqref{11} of Proposition~\ref{P} and the identity
%
%e90 ###
%e88 #&#
\begin{equation} \label{e_i_ddm}
\frac{1}{s} \frac{d}{d \sigma} = s \frac{d}{dm},
\end{equation}
which is a direct consequence of \eqref{36}. Let us consider the
second term. The identity
\[
\biggl ( \frac{1}{s} \frac{d}{d \sigma} \biggr)^2 \overset{\mathclap{\eqref
{e_i_ddm}}}{=}  \biggl( s \frac{d}{dm} \biggr) \biggl ( s \frac{d}{dm}
 \biggr) = s^2 \frac{d^2}{dm^2} + s \biggl ( \frac{d}{dm} s  \biggr)
\frac{d}{dm} ,
\]
which we rewrite as
\[
s^2 \frac{d^2}{dm^2} = \biggl ( \frac{1}{s} \frac{d}{d \sigma}
\biggr)^2 - \biggl ( \frac{d}{dm} s  \biggr)  \frac{1}{s} \frac{d}{d\sigma},
\]
yields
\[
\frac{d^2}{dm^2} g_{K, \sigma} (0) = \frac{1}{s^2} \biggl ( \biggl (
\frac{1}{s} \frac{d}{d \sigma} \biggr)^2 g_{K, \sigma} (0) -
\biggl(\frac{d}{dm} s  \biggr) \frac{1}{s} \frac{d}{d\sigma} g_{K, \sigma
}(0)  \biggr).
\]
Now, estimates \eqref{11} and \eqref{12} of Proposition~\ref{P} and
Lemma~\ref{e_aux_est_d_dm_s} yield the desired estimate \eqref
{e_d^2_dm_log_g}.
\end{pf*}

%apA #&#
\begin{appendix}
\section*{Appendix: Standard criteria for the SG and the LSI}\label
{s_basic_facts_about_LSI}
\setcounter{theorem}{0}

In this section we quote some standard criteria for the SG and the LSI.
For a general introduction to the SG and the LSI we refer to \cite
{L,R,GZ}. Note that even if we only formulate the criteria on the level
of the LSI, they also hold on the level of the SG. The first one shows
that the LSI is compatible with products; cf., for example, \cite{GZ}, Theorem 4.4.

%thA.1 #&#
\begin{theorem}[(Tensorization principle)]\label{p_tensorization_principle}
Let $\mu_1$ and $\mu_2$ be probability measures on Euclidean spaces
$X_1$ and $X_2$, respectively. If $\mu_1$ and $\mu_2$ satisfy the LSI
with constant $\varrho_1$ and $\varrho_2$, respectively, then the
product measure $\mu_1 \otimes\mu_2$ satisfies the LSI with constant
$\min\{\varrho_1, \varrho_2 \}$.
\end{theorem}

The next criterion shows how the LSI constant behaves under
perturbations; cf. \cite{HS}, page~1184.

%thA.2 #&#
\begin{theorem}[(Holley--Stroock criterion)] \label{p_criterion_holley_stroock}
Let $\mu$ be a probability measure on the Euclidean space $X$, and let
$\delta\psi\dvtx  X \to\R$ be a bounded function. Let the probability
measure $\tilde\mu$ be defined as
\[
\tilde\mu(dx) = \frac{1}{Z} \exp ( - \delta\psi(x) )
  \mu(dx).
\]
If $\mu$ satisfies the LSI with constant $\varrho$, then $\tilde\mu
$ satisfies the LSI with constant
\[
\tilde\varrho= \varrho \exp \bigl( -  ( \sup\delta\psi-
\inf\delta\psi ) \bigr).
\]
\end{theorem}

Because of its perturbative nature, the Holley--Stroock criterion is
not well adapted for high dimensions. For the proof of the last
statement, we refer the reader to \cite{L}, Lemma 1.2. Now, we state
the Bakry--\'Emery criterion, which connects the convexity of the
Hamiltonian to the LSI constant; cf. \cite{BE}, Proposition~3 and Corollary~2, or \cite{L}, Corollary~1.6.\vadjust{\goodbreak}
%
%thA.3 #&#
\begin{theorem}[(Bakry--\'Emery criterion)] \label{p_criterion_bakry_emery}
Let $d\mu:= Z^{-1} \exp(-H(x))\,dx$ be a probability measure on a
Euclidean spaces $X$. If there is a constant $\varrho>0$ such that in
the sense of quadratic forms
\[
\Hess H(x) \geq\varrho
\]
uniformly in $x \in X$, then $\mu$ satisfies the LSI with constant
$\varrho$.
\end{theorem}

A proof using semi-group methods can be found in \cite{L}, Corollary~1.6. There is also a heuristic interpretation of the
Bakry--\'Emery criterion on a formal Riemannian structure on the space
of probability measures; cf. \cite{OV}, Section 3.
\end{appendix}

\section*{Acknowledgment}
G. Menz and F. Otto thank Franck Barthe, Michel Ledoux and Cedric
Villani for discussions on this subject.

%suskaldyti doi

% imsref loaded by smiklovaite, 2012-03-07 13:16:36
%

\printaddresses

\end{document}